\documentclass[12pt]{article}

\usepackage[english]{babel}
\usepackage[T1]{fontenc}
\usepackage{amsmath}
\usepackage{amsfonts}
\usepackage{amssymb}
\usepackage[utf8]{inputenc}
\usepackage{graphicx,color,url}
\usepackage[margin=2.4cm]{geometry}
\usepackage{booktabs}
\usepackage{caption}
\usepackage{tikz,pgfplots}
        \pgfplotsset{compat = 1.3}
        \pgfplotsset{minor grid style={dotted}} \pgfplotsset{major grid
        style={dashed}}

        \pgfplotsset{every x tick label/.append style={font=\footnotesize,
        yshift=0.25ex}}

        \pgfplotsset{every y tick label/.append
        style={font=\footnotesize, xshift=0.25ex}}
\usetikzlibrary{shapes.arrows}

\newtheorem{rem}{Remark}

\def\mca{\mathcal{A}}
\def\mcq{\mathcal{Q}}

\def\pmatrix{\left(\begin{array}}
\def\endpmatrix{\end{array}\right)}
\def\dd{\mathrm{d}}

\def\sech{\mathrm{sech}}

\def\sech{\mathrm{sech}}
\def\p{\partial}
\def\mca{\mathcal{A}}
\def\mcq{\mathcal{Q}}
\def\wmca{\widetilde{\mathcal{A}}}
\def\wmcq{\widetilde{\mathcal{Q}}}

\begin{document}

\title{Numerical preservation of multiple local conservation laws}

\author{G.\,Frasca-Caccia\,\quad P.\,E.\,Hydon\,\\[.5cm]
\small
School of Mathematics, Statistics and Actuarial Science\\
\small
University of Kent, Canterbury, CT2 7NZ\\
\small}

\maketitle

\begin{abstract} 
There are several well-established approaches to constructing finite difference schemes that preserve global invariants of a given partial differential equation. However, few of these methods preserve more than one conservation law {\em locally}. A recently-introduced strategy uses symbolic algebra to construct finite difference schemes that preserve several local conservation laws of a given scalar PDE in Kovalevskaya form. In this paper, we adapt the new strategy to PDEs that are not in Kovalevskaya form and to systems of PDEs. The Benjamin--Bona--Mahony equation and a system equivalent to the nonlinear Schr\"odinger equation are used as benchmarks, showing that the strategy yields conservative schemes which are robust and highly accurate compared to others in the literature.


\textbf{AMS class:} \textit{65M06, 37K05, 39A14}
\end{abstract}

\section{Introduction}

In the numerical treatment of partial differential equations (PDEs), the benefits of preserving global integral invariants are well-known (see \cite{MatsuoFur,S-SV86,Heitz,KoideFur,BFI15,BarlettiNLS,DuranLM,DuranSS,
defrutosSS,DuranSS98,Celledoni,Dahlby,McLaren}). Such invariants are obtained by integrating (local) conservation laws of the PDE over a spatial domain and applying given boundary conditions. However, even for PDEs in conservation form, very few invariant-preserving finite difference schemes preserve more than one conservation law locally. The main exceptions are Hamiltonian PDEs, for which a fully discrete local conservation law for the Hamiltonian \cite{McLQuisp} can be obtained by first creating a spatial discretization that preserves a semidiscrete conservation law for the Hamiltonian, then using a discrete gradient method \cite{Gonz,Cell2,Robi,QuispTurn,Celledoni,Dahlby,McLaren} on the resulting system of ordinary differential equations.

To preserve local features of any PDE, a scheme must satisfy much stronger constraints than are needed to preserve the corresponding global properties. Conservation laws hold throughout the domain; they apply to the set of all solutions, so are independent of initial and boundary conditions. Moreover, conservation laws correspond to topological properties of the solution manifold: each equivalence class of conservation laws is a cohomology class of the variational bicomplex restricted to this manifold (see \cite{An1992,Vi1984}). There is a strong theme in geometric integration that numerical approximations should respect topological properties of a given PDE as far as possible, so it is worthwhile trying to develop schemes that preserve discrete analogues of multiple conservation laws.

A new approach that uses symbolic algebra to develop bespoke finite difference schemes that preserve multiple local conservation laws of a (not necessarily Hamiltonian) PDE was introduced in \cite{Grant,GrantHydon}. Initially, the complexity of the symbolic calculations made this approach impractical for all but the simplest PDEs. However, a strategy introduced recently in \cite{FCHydon} has overcome this difficulty, making it possible to create bespoke schemes for a given PDE within a few minutes. This strategy has been used to create robust, highly accurate conservative schemes for the KdV and modified KdV equations, and a (non-Hamiltonian) nonlinear heat equation \cite{FC18,FCHydon,FCHmkdv}. All of these PDEs are in Kovalevskaya form.

The main aim of this paper is to show that this strategy can also be used for PDEs that are not in Kovalevskaya form and systems of PDEs. As benchmark examples, we will create conservative schemes for the Benjamin--Bona--Mahony (BBM) equation (or Regularized Long Wave equation) and for a system of two real PDEs that is equivalent to the Nonlinear Schr\"odinger (NLS) equation. Both equations have been well-studied; each is Hamiltonian, with a wide range of schemes that preserve some of their geometric properties. 

Multisymplectic methods, which preserve a local conservation law for symplecticity \cite{Bridges,BridgesReich} were applied to the BBM equation in \cite{LiSun,QinSun} and to the NLS equation in \cite{Chen2000,IslasSchober,SunQin,ChenQin}. In general, multisymplectic methods do not preserve further conservation laws or global invariants. The benefits of using (global) invariant-preserving schemes for the BBM and the NLS equation were studied in \cite{DuranLM} and \cite{DuranSS}, respectively. Several such schemes for the BBM equation were proposed in \cite{KoideFur} and for the NLS equation in \cite{MatsuoFur,S-SV86,BarlettiNLS,Heitz}. Using the inverse scattering transform, Ablowitz and Ladik \cite{AblowitzLadik} found one of the most important semi-discrete models of the NLS equation; it is completely integrable. Its symplectic structure is noncanonical, so standardization is needed before symplectic integrators can be used \cite{Schober,Tang}.

This paper is organized as follows: in Section~\ref{method} we present the basic strategy, adapted to systems of PDEs that are not necessarily in Kovalevskaya form. In Section~\ref{BBMsec} we use this strategy to find finite difference schemes that preserve two conservation laws of the BBM equation. In Section~\ref{NLSsec} the strategy is applied to the NLS equation. Time discretizations of the Ablowitz--Ladik model that preserve two local conservation laws of the NLS are also derived. A range of numerical tests show the effectiveness of the proposed schemes and compare them with known geometric integrators. We conclude with some remarks in Section~\ref{Conclusions}.

\section{Constructing schemes that preserve conservation laws}\label{method}

A system of PDEs for $\mathbf{u}=(u^1(x,t),\dots,u^q(x,t))$ may be written in the form
\begin{equation}\label{PDE}
\mathcal{A}(x,t,[\mathbf{u}])=\mathbf{0},
\end{equation}
where $\mathcal{A}(x,t,[\mathbf{u}])$ is a row vector and $[\mathbf{u}]$ denotes $\mathbf{u}$ and finitely many of its derivatives. Similarly, square brackets around any differentiable expression are used to denote the expression and a finite number of its derivatives. We assume throughout that the system of PDEs is totally nondegenerate, avoiding pathological exceptions (see \cite{Olverbook} for details).

A (local) conservation law of (\ref{PDE}) is a total divergence,
$$\text{Div}\,\mathbf{F}\equiv D_x\{F(x,t,[\mathbf{u}])\}+D_t\{G(x,t,[\mathbf{u}])\},$$
that vanishes on all solutions of (\ref{PDE}), so that
\begin{equation}\label{contCLaw}
\text{Div}\,\mathbf{F}=0\quad\text{when}\quad[\mathcal{A}=\mathbf{0}].
\end{equation}
Here
\begin{align*}
D_x&=\frac{\p}{\p x}+u^\alpha_x\frac{\p}{\p u^\alpha}+u^\alpha_{xt}\frac{\p}{\p u^\alpha_t}+u^\alpha_{xx}\frac{\p}{\p u^\alpha_x}+\cdots,\\
D_t&=\frac{\p}{\p t}+u^\alpha_t\frac{\p}{\p u^\alpha}+u^\alpha_{tt}\frac{\p}{\p u^\alpha_t}+u^\alpha_{xt}\frac{\p}{\p u^\alpha_x}+\cdots
\end{align*}
are the total derivatives with respect to $x$ and $t$, respectively, and the functions $F$ and $G$ are the flux and the density of the conservation law, respectively. A conservation law is in \textit{characteristic form} if
\begin{equation}\label{charclaw}
\text{Div}\,\mathbf{F}=\mathcal{A}\mathcal{Q},
\end{equation}
for some column vector $\mathcal Q$, which is called the characteristic. (Every conservation law for a totally nondegenerate system is equivalent to one in characteristic form.)

The kernel of the Euler operator $\mathcal{E}$, whose $\alpha$-th component is
\[\mathcal{E}_\alpha=\sum_{i,j}(-D_x)^i(-D_t)^j\frac{\partial}{\partial u_{x^it^j}^\alpha},\quad\text{where}\quad u_{x^it^j}^\alpha=D_x^iD_t^j(u^\alpha),\]
is the vector space of total divergences. Consequently, if
\[\mathcal{E}(\mathcal{AQ})=\mathbf{0},\]
there exists $\mathbf{F}$ such that $\mathcal{AQ}=\text{Div}\,\mathbf{F}$ is a conservation law.

For simplicity, we shall discretize the PDE (\ref{PDE}) and its conservation laws (\ref{contCLaw}) on a uniform lattice. Relative to a generic lattice point $\mathbf{n}=(m,n)$, the grid points are 
\begin{equation*}
 x_i=x(m+i)=x(m)+i\Delta x,\qquad t_j=t(n+j)=t(n)+j\Delta t.
\end{equation*}
The approximated values of $\mathbf{u}$ at these points are $\mathbf{u}_{i,j}$ with $\alpha$-th component
\begin{equation*}
u^\alpha_{i,j}\approx u^\alpha(x_i,t_j),\qquad i,j\in\mathbb{Z}.
\end{equation*}
The forward shift operators in space and time are defined by
$$S_m:f(x_i,t_j)\mapsto f(x_{i+1},t_j),\qquad S_n:f(x_i,t_j)\mapsto f(x_{i},t_{j+1}),$$
for any function $f$ that is defined on the grid. The forward difference operators $D_m$, $D_n$ and the forward average operators $\mu_m$, $\mu_n$ are
\begin{eqnarray*}
D_m=\tfrac{1}{\Delta x}(S_m-I),\quad D_n=\tfrac{1}{\Delta t}(S_n-I),\quad
\mu_m=\tfrac{1}{2}(S_m+I),\quad \mu_n=\tfrac{1}{2}(S_n+I),
\end{eqnarray*}
where $I$ is the identity operator.

Discretizing (\ref{PDE}) by a suitable finite difference approximation yields a system of partial difference equations (P$\Delta$Es),
\begin{equation*}
\widetilde{\mathcal{A}}(m,n,[\mathbf{u}])=\mathbf{0}.
\end{equation*}
Here $[\mathbf{u}]$ denotes $\mathbf{u}_{0,0}$ and a finite number of its shifts; more generally, square brackets around a discrete expression denote the expression and finitely many of its shifts. Here and henceforth tildes represent discretizations of the corresponding continuous terms. (See \cite{Hydonbook} for a comprehensive introduction to difference equations and their conservation laws.)

We seek schemes having the following discrete analogue of each preserved conservation law:
\begin{equation}\label{discCLaw}
\text{Div}\,\widetilde{\mathbf{F}}\equiv D_m\{\widetilde{F}(m,n,[\mathbf{u}])\}+D_n\{\widetilde{G}(m,n,[\mathbf{u}])\},
\end{equation}
such that $$\text{Div}\,\widetilde{\mathbf{F}}=0 \quad \text{when}\quad [\widetilde{\mathcal{A}}=\mathbf{0}].$$
The functions $\widetilde{F}$ and $\widetilde{G}$ are called the discrete flux and the discrete density of the conservation law (\ref{discCLaw}), respectively.
A discrete conservation law is said to be in characteristic form if there exists $\widetilde{\mathcal{Q}}$, called the characteristic, such that
\[D_m\{\widetilde{F}(m,n,[\mathbf{u}])\}+D_n\{\widetilde{G}(m,n,[\mathbf{u}])\}=\widetilde{\mathcal{A}}\widetilde{\mathcal{Q}}(m,n,[\mathbf{u}]).\]
Just as for PDEs, every conservation law for a totally nondegenerate system of difference equations is equivalent to one in characteristic form \cite{Hydonbook}.

The linear and quadratic terms in $\mca$ and $\mcq$ are, respectively, approximated by
\begin{align}\label{linapp}
\frac{\partial^{r+s}}{\partial x^r\partial t^s}u^\alpha&\approx\frac{1}{\Delta x^r}\frac{1}{\Delta t^s}\sum_{i=A}^B\sum_{j=0}^1\alpha_{i,j}u_{i,j}^\alpha,\\\label{quadapp}
\frac{\partial^{p+q}u^\alpha}{\partial x^p\partial t^q}\frac{\partial^{r+s}u^\alpha}{\partial x^r\partial t^s}&\approx\frac{1}{\Delta x^{p+r}\Delta t^{q+s}}\sum_{i=A}^B\left(\sum_{k=i}^B\sum_{j=0}^1\beta_{i,j,k}u^\alpha_{i,j}u^\alpha_{k,j}+\sum_{k=A}^B\gamma_{i,k}u^\alpha_{i,0}u^\alpha_{k,1}\right).
\end{align}
Here the coefficients $\alpha_{i,j}$, $\beta_{i,j,k}$ and $\gamma_{i,k}$ (which depend on $p,q,r,s,\alpha$) are chosen to ensure that the Taylor expansion of each approximated term about the centre of the stencil is accurate to whatever order, ${\rho}$, is desired. Provided that the stencil is sufficiently large, the approximations of $\mca$ and $\mcq$ will include some free parameters. (For a given stencil, the number of free parameters decreases rapidly with the order of the approximation.) Terms that are differential polynomials of third and higher degree are approximated similarly.

The following result is the key to obtaining conservative schemes (see \cite{Kuperschmidt}, with generalizations in \cite{HydonMans}).

\begin{rem}\label{eulerremark}
The kernel of the difference Euler operator $\mathsf{E}$, whose $\alpha$-th component is
\[\mathsf{E}_\alpha=\sum_{i,j}S_m^{-i}S_n^{-j}\frac{\partial}{\partial u^\alpha_{ij}}\,,\]
coincides with the space of difference divergences (\ref{discCLaw}).
\end{rem}

Given a PDE (\ref{PDE}) with a conservation law in characteristic form (\ref{charclaw}) that one wishes to preserve, the aim is to seek approximations $\widetilde{\mcq}$ and $\widetilde{\mca}$ such that
\[\mathsf{E}(\mathcal{\widetilde A\widetilde Q})=\mathbf{0}.\]
By Remark~\ref{eulerremark}, there exist a discrete flux $\widetilde{F}$ and density $\widetilde{G}$ such that $\mathcal{\widetilde A\widetilde Q}=D_m \widetilde{F}+D_n\widetilde{G}$ is a difference conservation law that approximates the corresponding continuous one.

Our strategy for obtaining bespoke conservative finite difference schemes is as follows.
\begin{enumerate}
	\item Select the desired accuracy, $\rho$, and the conservation laws to be preserved, labelling their characteristics $\mcq_1, \mcq_2,\dots$. We shall require second-order accuracy, to reduce the number of free parameters to a manageable level while retaining enough freedom to preserve multiple conservation laws.
	\item Choose a rectangular stencil that is large enough to support ${\rho}^{\mathrm{th}}$-order approximations of $\mca$ and all $\mcq_\ell$.
	\item Determine the most general ${\rho}^{\mathrm{th}}$-order finite difference approximations to $\wmca$ and $\wmcq_1$ on the stencil.
	\item
	Reduce the number of free parameters remaining by making some key terms as compact as possible; typically, these include highest derivatives and highest-order nonlinear terms.
	\item Use symbolic algebra to solve for the values of the free parameters that satisfy
	\begin{equation}\label{Eulcond}
	\mathsf{E}(\wmca\wmcq_1)=0.
	\end{equation}
	For each solution (which may depend on free parameters), the discrete flux $\widetilde{F}_1$ and density $\widetilde{G}_1$ can be reconstructed from the characteristic (see \cite{GrantHydon,Hydon2001}).
	\item Iterate Steps 3 onwards (with $\mcq_\ell$ replacing $\mcq_1$, etc.) to preserve further conservation laws. If $\mathsf{E}(\wmca\wmcq_\ell)=0$ has no solutions for some $\ell$, the corresponding conservation law cannot be preserved on the chosen stencil without violating an earlier conservation law.
\end{enumerate}

The restriction to second-order approximations and the compactness conditions were introduced in \cite{FCHydon} to obtain solutions of (\ref{Eulcond}) by means of a fast symbolic computation. On the stencils we consider below, this takes around one minute on a fast laptop. By contrast, if one imposes no constraints other than convergence, the symbolic computation typically takes several days; this is too long to be practical.

\begin{rem}\label{rem2ord}
	Obtaining second-order accurate approximations of the conservation laws at the centre $(x, t)$ of the stencil is equivalent to finding second-order accurate approximations of the corresponding densities and fluxes at the points $(x, t-\Delta t/2)$ and $(x-\Delta x/2, t)$ respectively.
\end{rem}

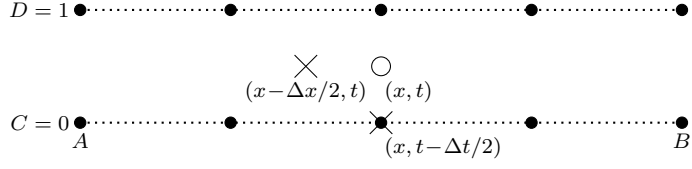
\begin{figure}
\center{
\begin{tikzpicture}
\draw[fill] (0,3) circle [radius=0.075];
\draw[fill] (0,1.5) circle [radius=0.075];
\draw[fill] (2,1.5) circle [radius=0.075];
\draw[fill] (4,1.5) circle [radius=0.075];
\draw[fill] (6,1.5) circle [radius=0.075];
\draw[fill] (8,1.5) circle [radius=0.075];
\draw[fill] (2,3) circle [radius=0.075];
\draw[fill] (4,3) circle [radius=0.075];
\draw[fill] (6,3) circle [radius=0.075];
\draw[fill] (8,3) circle [radius=0.075];
\draw[thick,dotted] (0,1.5)--(8,1.5);
\draw[thick,dotted] (0,3)--(8,3);
\node [below right] at (3.9,1.45) {\scriptsize{$(x,t\!-\!{\Delta t/2})$}};
\node [below right] at (3.9,2.2) {\scriptsize{$(x,t)$}};
\node [below] at (3,2.2) {\scriptsize{$(x\!-\!{\Delta x/2},t)$}};
\draw (4,2.25) circle [radius=0.125];
\draw (3.85,1.35)--(4.15,1.65);
\draw (3.85,1.65)--(4.15,1.35);
\draw (2.85,2.10)--(3.15,2.4);
\draw (2.85,2.4)--(3.15,2.1);
\node [below] at (0,1.5) {\scriptsize{$A$}};
\node [below] at (8,1.5) {\scriptsize{$B$}};
\node [left] at (0,1.5) {\scriptsize{$C=0$}};
\node [left] at (0,3) {\scriptsize{$D=1$}};
\end{tikzpicture}
\caption{Example of a rectangular stencil for one-step schemes. PDEs and conservation laws are preserved to second order at the central point $(x, t)$; densities and fluxes are second-order at $(x, t-\Delta t/2) $ and $(x-\Delta x/2, t)$, respectively.}
\label{stencil}
}
\end{figure}

In the next two sections, we use our strategy to develop conservative finite difference schemes for the BBM equation and for the system of two PDEs given by the real formulation of the NLS equation. As these are evolution equations, we will consider only one-step schemes defined on the stencil in Figure~\ref{stencil}, where $B-A$ is suitably large. However, the strategy applies equally to multistep schemes.

\section{The Benjamin--Bona--Mahony equation}\label{BBMsec}
In this section we consider the BBM equation,
\begin{equation}\label{BBM}
\mca\equiv u_t-uu_x-u_{xxt}=0, \quad (x,t)\in\Omega\equiv[a,b]\times[0,\infty),
\end{equation}
which (up to equivalence) has just three independent conservation laws \cite{Olver,DuzhinTsujishita}:
\begin{align}\label{BBMCL1}
& D_x(F_1)+D_t(G_1) = D_x\left(-\tfrac{1}2 u^2 - u_{xt}\right)+D_t(u)=0,\\\label{BBMCL2}
& D_x(F_2)+D_t(G_2) = D_x\left(-\tfrac{1}3u^3-uu_{xt}\right)+D_t\left(\tfrac{1}2 u^2+\tfrac{1}2u_x^2\right)=0,\\\label{BBMCL3}
& D_x(F_3)+D_t(G_3) = D_x\left(u_t^2-\tfrac{1}4u^4-u^2u_{xt}-u_{xt}^2\right)+D_t\left(\tfrac{1}3u^3\right)=0.
\end{align}
These can be written in characteristic form with characteristics
\[
Q_1=1,\qquad Q_2=u,\qquad Q_3=u^2+2u_{xt},
\]
respectively. When (\ref{BBM}) is coupled with suitable (e.g. periodic) boundary conditions, integrating (\ref{BBMCL1})--(\ref{BBMCL3}) over the spatial domain gives the invariants:
\begin{equation}\label{glinvbbm}
\int G_1\,\dd x=\int u\,\dd x,\quad \int G_2\,\dd x=\int \tfrac{1}{2}(u^2+u_x^2)\,\dd x, \quad \int G_3\,\dd x=\int \tfrac{1}{3}u^3\, \dd x,
\end{equation}
which are usually referred to as the ``mass'', ``momentum'' and ``energy'', respectively \cite{KoideFur}.

Equation (\ref{BBM}) admits the Hamiltonian formulation,
\begin{equation*}
u_t=\left(1-\frac{\partial^2}{\partial x^2}\right)^{-1}\frac{\partial}{\partial x}\left(\frac{\delta}{\delta u}\mathcal{H}\right),
\end{equation*}
where ${\delta}/{\delta u}$ is the variational derivative and the Hamiltonian functional $\mathcal{H}$ is given by
\begin{equation*}
\mathcal{H}=\int\tfrac{1}6 u^3\,\dd x=\int \tfrac{1}2G_3\,\dd x.
\end{equation*}
Given suitable boundary conditions, the local conservation law (\ref{BBMCL3}) implies the conservation of the Hamiltonian functional.
\subsection{Conservative methods for the BBM equation}\label{BBMmeth}
This section describes how our strategy produces numerical schemes for the BBM equation that preserve two local conservation laws. All of these methods are of the form
\begin{equation}\label{genmet}
\widetilde{\mca}=D_m\widetilde{F}_1+D_n\widetilde{G}_1=0,
\end{equation}
so the mass conservation law (\ref{BBMCL1}) is preserved. Some schemes depend on free parameters, all of which are $\mathcal{O}(\Delta x^2,\Delta t^2)$. As we impose only second-order accuracy, one can find values of the parameters that reduce the local truncation error. However, no choice of the parameters yields any higher-order scheme; indeed, their optimal values depend on the particular problem. We use three different stencils but none of the resulting schemes preserves all three conservation laws of the BBM equation.

\subsubsection*{6-point schemes}
The most compact stencil for the BBM equation has 6 points. We choose $A=-1$, $B=1$ in Figure~\ref{stencil} and seek second-order approximations of characteristics, densities and fluxes at $(0,1/2)$, $(0,0)$ and $(-1/2,1/2)$, respectively.

\quad\\
\textbf{Energy-conserving methods.}\\
The symbolic calculation that finds 6-point energy-conserving schemes is simplified here by specifying the approximations of the highest derivatives in ${F_1}$ and ${\mathcal{Q}_3}$ respectively to be
\begin{align}\label{app6st}
\widetilde{u_{x,t}}&= D_mD_nu_{-1,0},\\\nonumber
\widetilde{u_{x,t}}&=D_mD_n\mu_m u_{-1,0}.
\end{align}
Approximations of linear and quadratic terms in $F_1$, $G_1$ and $\mathcal{Q}_3$ are of the form (\ref{linapp}) and (\ref{quadapp}), subject to the requirement of second-order accuracy. Energy conservation is obtained by solving
\begin{equation*}
\mathsf{E}(\widetilde{\mca}\widetilde{\mcq}_3)\equiv 0,
\end{equation*}
which determines all of the remaining coefficients in $\widetilde{F}_1$, $\widetilde{G}_1$ and $\widetilde{\mathcal{Q}}_3$, giving scheme $\mbox{EC}_6$ in the form (\ref{genmet}), with
\begin{equation*}
\widetilde{F}_1=-\tfrac{1}{6}\left(\left(\mu_m u_{-1,1}\right)^2+\left(\mu_m u_{-1,0}\right)^2+\left(\mu_m u_{-1,1}\right)\left(\mu_m u_{-1,0}\right)\right)-D_mD_nu_{-1,0},\qquad\widetilde{G}_1=u_{0,0}.
\end{equation*}
The scheme $\mbox{EC}_6$ preserves the following discrete version of the conservation law (\ref{BBMCL3}):
\begin{align*}
\widetilde{\mca}\widetilde{\mcq}_3=& D_m\widetilde{F}_3+ D_n\widetilde{G}_3=0,
\end{align*}
with
\begin{align}\nonumber
\widetilde{F}_3=&-\widetilde{F}_1^2+(D_n u_{-1,0})D_n u_{0,0}-\tfrac{1}{3}\Delta x^2(\mu_m\mu_nu_{-1,0})\Theta[u],\\\nonumber
\widetilde{G}_3=&\,\tfrac{1}{3}{u_{0,0}\mu_m\left(\left(\mu_mu_{-1,0}\right)^2\right)},\qquad \widetilde{\mcq}_3=-2\mu_m\widetilde{F}_1, 
\end{align}
and
\begin{equation}\label{ThetaBBM}
\Theta[u]=\tfrac{1}2\left\{(\mu_m\mu_nu_{-1,0})D_mD_nu_{-1,0}-(D_m\mu_nu_{-1,0})D_n\mu_mu_{-1,0}\right\}.
\end{equation}
The last term in $\widetilde{F}_3$ does not correspond to an expression in the continuous flux; it vanishes as the spatial stepsize tends to zero.\\
\quad\\
\textbf{Momentum-conserving methods.}\\
Momentum-conserving methods are obtained by specifying the approximation of the highest derivative in $F_1$ as in (\ref{app6st}) and using the compact approximation
\[
\widetilde{\mcq}_2=\mu_n u_{0,0}.
\]
The remaining terms in $F_1$ and $G_1$ are approximated by (\ref{linapp}) and (\ref{quadapp}). Solving
\begin{equation*}
\mathsf{E}(\widetilde{\mca}\widetilde{\mcq}_2)\equiv 0
\end{equation*}
gives an one-parameter family of mass and momentum-conserving methods $\widetilde{\mathcal{A}}(\lambda)$ with
\begin{align*}
\widetilde{F}_1=-\tfrac{1}{6}\left(\left(\mu_n u_{-1,0}\right)^2+\left(\mu_n u_{0,0}\right)^2+\left(\mu_n u_{-1,0}\right)\left(\mu_n u_{0,0}\right)\right)+(\lambda-1)D_mD_nu_{-1,0},\quad\widetilde{G}_1=u_{0,0},
\end{align*}
and $\lambda=\mathcal{O}(\Delta x^2,\Delta t^2)$.
For any value of $\lambda$, these methods preserve the discrete momentum conservation law
$$\widetilde{\mca}\widetilde{\mcq}_2= D_m\widetilde{F}_2+ D_n\widetilde{G}_2=0,$$
with
\begin{align*}
\widetilde{F}_2=&\,-\tfrac{1}{3}(\mu_nu_{-1,0})(\mu_nu_{0,0})\mu_m\mu_nu_{-1,0}-(\mu_m\mu_nu_{-1,0})D_mD_nu_{-1,0}+\lambda\Theta[u],\\
\widetilde{G}_2=&\,\tfrac{1}{2}\left\{u_{0,0}^2+\mu_m\left((D_m u_{-1,0})^2\right)+\lambda u_{0,0}D_m^2u_{-1,0}\right\},
\end{align*}
and $\Theta[u]$ as given in (\ref{ThetaBBM}).  In the numerical tests, we use the notation\[\mbox{MC}_6(\alpha)=\widetilde{\mca}(\alpha\Delta_{\text{max}}^2),\qquad \Delta_{\text{max}}=\max{(\Delta x,\Delta t)}.\]

\subsubsection*{8-point schemes}
We now develop schemes that preserve two conservation laws of the BBM equation, defined on the 8-point stencil with $A=-2$, $B=1$ in Figure~\ref{stencil}. Approximations of characteristics, densities and fluxes are second-order accurate at $(-1/2,1/2)$, $(-1/2,0)$ and $(-1,1/2)$, respectively.

\quad\\
\textbf{Energy-conserving methods.}\\
Mass and energy-conserving schemes are obtained here by considering the approximation (\ref{app6st}) of $u_{x,t}$ in $F_1$ and by considering only compact approximations of the quadratic term in ${Q_3}$.

All the approximations of the remaining terms in $F_1$, $G_1$ and $\mcq_3$ are in the form (\ref{linapp}) or (\ref{quadapp}). The condition
\begin{equation*}
\mathsf{E}(\widetilde{\mca}\widetilde{\mcq}_3)\equiv 0,
\end{equation*}
can be solved with a fast symbolic computation, yielding scheme $\mbox{EC}_{8}$ in the form (\ref{genmet}) with flux and density
\begin{equation*}
\widetilde{F}_1=\mu_m\varphi_{-2,0},\quad\widetilde{G}_1=\mu_m u_{-1,0},
\end{equation*}
and
$$\varphi_{-2,0}=-\tfrac{1}{6}\left(\left(\mu_m u_{-2,1}\right)^2+\left(\mu_m u_{-2,0}\right)^2+\left(\mu_m u_{-2,1}\right)\left(\mu_m u_{-2,0}\right)\right)-D_mD_nu_{-2,0}.$$
The scheme $\mbox{EC}_{8}$, preserves the following discrete version of the conservation law (\ref{BBMCL3}):
\begin{align*}
\widetilde{\mca}\widetilde{\mcq}_3=& D_m\widetilde{F}_3+ D_n\widetilde{G}_3=0,
\end{align*}
with
\begin{align*}
\widetilde{\mcq}_3=\,-2\varphi_{-1,0},\quad \widetilde{F}_3=\,-\varphi_{-2,0}\varphi_{-1,0}+(D_nu_{-1,0})^2,\quad \widetilde{G}_3=\,\tfrac{1}{3}(\mu_mu_{-1,0})^3.
\end{align*}

Note that $\mbox{EC}_8$ amounts to averaging $\mbox{EC}_6$ in space, previously introduced. No other $8$-point energy-preserving schemes can be found by making a different compactness assumption.

\quad\\
\textbf{Momentum-conserving methods.}\\
In order to reduce the complexity of the symbolic computations yielding momentum-conserving schemes on the 8-point stencil, we use the approximations in (\ref{app6st}) and the most compact approximations of $\widetilde{\mcq}_2$,
\[
\widetilde{\mathcal Q}_2=\mu_m\mu_nu_{-1,0},\\
\]
whereas all the other terms in $F_1$ and $G_1$ are approximated as in (\ref{linapp}) and (\ref{quadapp}). The solution of $$\mathsf{E}(\widetilde\mca\widetilde{\mcq}_2)\equiv 0,$$ gives two-parameter family of momentum-conserving methods $\widetilde\mca(\lambda,\nu)$ in the form (\ref{genmet}) with 
\begin{align*}\nonumber
\widetilde{F}_1=&\,-\tfrac{1}6(\mu_nu_{-1,0})\mu_n(u_{-2,0}+u_{-1,0}+u_{0,0})+(\lambda-1)D_nD_m\mu_mu_{-2,0}\\\nonumber
&\,+\nu\left\{\left(D_m\mu_nu_{-2,0}\right)D_m\mu_nu_{-1,0}+\left(D_m^2\mu_nu_{-2,0}\right)\mu_n(u_{-2,0}+u_{0,0})\right\},\\
\widetilde{G}_1=&\,\mu_mu_{-1,0},\\
\end{align*}
where $\lambda=\mathcal{O}(\Delta x^2,\Delta t^2)$ and $\nu=\mathcal{O}(\Delta x^2,\Delta t^2)$.
Each of these schemes preserves
\begin{align*}
\widetilde{\mca}\widetilde{\mcq}_2=& D_m\widetilde{F}_2+ D_n\widetilde{G}_2=0,
\end{align*}
with
\begin{align*}
\widetilde{F}_2=&\,-\tfrac{1}3(\mu_nu_{-1,0})(\mu_m\mu_nu_{-2,0})\mu_m\mu_nu_{-1,0}-\left(\mu_m^2\mu_n u_{-2,0}\right)D_mD_n\mu_mu_{-2,0}\\\nonumber
&\, +\lambda\left\{(\mu_m^2\mu_nu_{-2,0})D_mD_n\mu_mu_{-2,0}-(D_n\mu_m^2u_{-2,0})D_m\mu_m\mu_nu_{-2,0}\right\}\\
&\, -\nu(\mu_m\mu_nu_{-2,0})(\mu_m\mu_nu_{-1,0})D_m^2\mu_nu_{-2,0},\\\nonumber
\widetilde{G}_2=&\,\tfrac{1}2(\mu_mu_{-1,0})^2+\tfrac{1}2\mu_m\left(\left(D_m\mu_mu_{-2,0}\right)^2\right)+\lambda(\mu_mu_{-1,0})D_m^2\mu_mu_{-2,0}.
\end{align*}
We denote this two-parameter family of schemes by $$\mbox{MC}_8(\alpha,\beta)=\widetilde{\mca}(\alpha\Delta_{\text{max}} ^2,\beta\Delta_{\text{max}}^2),\qquad \Delta_{\text{max}}=\max{(\Delta x,\Delta t)}.$$

\subsubsection*{10-point energy-conserving schemes}
We now seek schemes that preserve local mass and energy on the 10-point stencil with $A=-2$, $B=2$ and second-order approximations of characteristics, densities and fluxes at $(0,1/2)$, $(0,0)$ and $(-1/2,1/2)$, respectively. In order to reduce the complexity of the symbolic calculation that solves $$\mathsf{E}(\widetilde{\mathcal{A}}\widetilde{\mathcal{Q}}_3)=0,$$ we use the following compact approximations of $G_1$ and the quadratic term in ${\mcq_3}$:
\begin{align*}
\widetilde{G}_1=u_{0,0},\qquad 
\widetilde{u^2}=\xi(u_{0,0}^2+u_{0,1}^2)+(1-2\xi)u_{0,0}u_{0,1};
\end{align*}
approximations of $F_1$ and the linear term in $\mathcal{Q}_3$ are of the form (\ref{linapp}) and (\ref{quadapp}).
This gives a family of schemes, $\widetilde{\mca}(\lambda)$, with
\begin{equation*}
\widetilde{F}_1=\mu_m\varphi_{-1,0},\quad \widetilde{G}_1= u_{0,0},
\end{equation*}
where
$$\varphi_{-1,0}=-\tfrac{1}6\left(u_{-1,0}^2+u_{-1,1}^2+u_{-1,0}u_{-1,1}\right)-D_mD_n\mu_m u_{-2,0}-\lambda D_m^2\mu_n u_{-2,0}$$
and $\lambda=\mathcal{O}(\Delta x^2,\Delta t^2).$
For any value of $\lambda,$ the schemes preserve
\begin{align*}
\widetilde{\mca}\widetilde{\mcq}_3=& D_m\widetilde{F}_3+ D_n\widetilde{G}_3=0,
\end{align*}
with
\[
\widetilde{\mcq}_3=-2\varphi_{0,0},\quad \widetilde{F}_3=-\varphi_{0,0}\varphi_{-1,0}+(D_nu_{0,0})D_nu_{-1,0}-2\lambda\Theta[u],\quad
\widetilde{G}_3=\tfrac{1}{3}u_{0,0}^3+\lambda u_{0,0}D_m^2u_{-1,0},
\]
where $\Theta[u]$ is given in (\ref{ThetaBBM}).
We use the notation
$$\mbox{EC}_{10}(\alpha)=\widetilde{\mca}(\alpha\Delta_{\text{max}}^2), \qquad \Delta_{\text{max}}=\max{(\Delta x,\Delta t)}.$$
The scheme $\mbox{EC}_{10}(0)$ amounts to the scheme derived by Koide and Furihata in \cite{KoideFur} using a Discrete Variational Derivative Method (DVDM).

\subsection{Numerical tests}

In this section we consider two benchmark problems in order to compare the schemes in Section~\ref{BBMmeth} with two known multisymplectic and mass-conserving schemes. These are the scheme proposed by Li and Sun in \cite{LiSun},
\begin{align*}
\text{LS}=D_n(\mu_mu_{-1,0})+D_m\left(-\tfrac{1}{2}(\mu_nu_{-1,0})^2-D_nD_m\mu_mu_{-2,0}\right)=0,
\end{align*}
and the multisymplectic Preissman box scheme \cite{Bridges,LiSun,Preissman,SunQin,QinSun,BridgesReich}, 
\begin{align*}
\mbox{PB}=D_n\left(\mu_m^3\mu_nu_{-2,-1}\right)+D_m\left(\mu_m\mu_n\left(-\tfrac{1}{2}(\mu_m\mu_nu_{-2,-1})^2-D_nD_mu_{-2,-1}\right)\right)=0.
\end{align*}
The scheme LS is defined on the 8-point stencil, whereas PB is a two-step method defined on a wider stencil consisting of 12 points.

For each of our numerical experiments, the computational time is similar for all of the schemes, so the main difference is the error in the solution at the final time $t=T$, evaluated as
\begin{equation}\label{solerr}
\left.\frac{\|u-u_{\mathrm{exact}}\|}{\|u_{\mathrm{exact}}\|}\right\vert_{t=T}.
\end{equation}

In all the experiments the BBM equation (\ref{BBM}) is defined on the domain $\Omega=[a,b]\times[0,T]$ equipped with periodic boundary conditions. 
Considering a grid with $M$ points in space and $N$ in time, we evaluate the error in the conservation laws by measuring the error in the discrete global invariants:
\begin{align}\label{erralpha}
\text{Err}_\alpha=&\,\Delta x\max_{j=1,\ldots,N}\left|\sum_{i=1}^M \left(\widetilde{G}_\alpha(x_i,t_j)-\widetilde{G}_\alpha(x_i,t_1)\right)\right|,\qquad \alpha=1,2,3.
\end{align}

Each discrete global invariant, $\sum_{i=1}^M\widetilde{G}_\alpha(x_i,t_j)$, is obtained by summing the corresponding conservation law over the spatial grid. The resulting approximation of the integrals in (\ref{glinvbbm}) can also be obtained by combining the composite trapezium quadrature rule and the second-order approximations of $G_\alpha$ obtained in Section~\ref{BBMsec}.

If $\widetilde{G}_2$ and/or $\widetilde{G}_3$ are not defined for a scheme, we evaluate the corresponding error as
\begin{align*}
\text{Err}_2=&\,\Delta x\max_{j=1,\ldots,N}\left|\tfrac{1}2\sum_{i=1}^M  \left(v_{i,j}^2+\mu_m\left((D_mv_{i-1,j})^2\right)-v_{i,0}^2-\mu_m\left((D_mv_{i-1,0})^2\right)\right)\right|,\\
\text{Err}_3=&\,\Delta x\max_{j=1,\ldots,N}\left|\tfrac{1}3\sum_{i=1}^M \left( {v_{i,j}^3}-{v_{i,0}^3}\right)\right|,
\end{align*}
with $v_{i,j}=u_{i,j}$ for schemes defined on the 6-point or 10-point stencil, or $v_{i,j}=\mu_mu_{i-1,j}$ otherwise. Subscripts denote shifts from $(x,t)=(a,0)$, so $u_{i,j}\simeq u(a+i\Delta x, j\Delta t)$.
We first consider (\ref{BBM}) on $\Omega=[-40,40]\times[0,T]$, with the initial condition obtained from the single solitary wave solution over $\mathbb{R}$:
\begin{equation}\label{ex1solbbm}
u_{\mathrm{exact}}(x,t)=3c\,\sech^2\!\left(\tfrac{1}2(x+ct-d)\right).
\end{equation}

It is useful to start by examining how errors accumulate over a long time, so we first solve the problem with $T=100$, $c=0.5$, $d=25$, stepsizes $\Delta x=0.05$ and $\Delta t=0.1$. The top-left plot in Figure~\ref{errorgrow} shows that the solution error (\ref{solerr}) of both MC$_8(0,0)$ and EC$_{10}(0)$ grows linearly with time; the solution errors for  EC$_6$, EC$_8$ and MC$_6(0)$ also grow linearly, each at a similar rate to that of MC$_8(0,0)$. By contrast, the solution error of a numerical method that does not preserve invariants generally accumulates quadratically (see, e.g., \cite{defrutosSS, DuranLM}). Among the five schemes considered, EC$_{10}(0)$ is the one that performs best in long time simulations, as the error accumulates at the slowest rate. The other plots in Figure~\ref{errorgrow} show the error in the conserved invariants at $t=0,0.1,\ldots,100$; this is due only to round-off errors. The error in the invariants conserved by MC$_8(0,0)$ and EC$_{10}(0)$ does not show any obvious drift due to accumulation of round-off errors; the same is true for EC$_8$. Linear numerical drift can be seen for EC$_6$, and MC$_6(0)$ has a similar drift.

\begin{figure}[t]
\begin{center}
	\definecolor{mycolor1}{rgb}{0.49020,0.18039,0.56078}%

\begin{tikzpicture}

\begin{axis}[%
width=2in,
height=1.5in,
at={(0.758in,0.481in)},
scale only axis,
xmin=0,
xmax=100,
xlabel style={font=\color{white!15!black}},
xlabel={$t$},
ymin=0,
tick label style={/pgf/number format/fixed},
scaled y ticks=false,
ymax=0.008,
ylabel style={font=\color{white!15!black}},
title={Solution error ($\times 10^{-3}$)},
ytick={0, 0.002, 0.004,  0.006,  0.008},
yticklabels={0, 2, 4, 6, 8},
axis background/.style={fill=white},
legend style={legend cell align=left, align=left, draw=white!15!black, at={(0.57,0.95)}}
]

\addplot [color=blue, line width=1.0pt, mark size=3.0pt, mark=o, mark options={solid, blue}]
  table[row sep=crcr]{%
5	0.000340794647322845\\
10	0.000682014803210115\\
15	0.0010235288125296\\
20	0.00136297108647519\\
25	0.00170134211237813\\
30	0.0020394370267989\\
35	0.00237748235873519\\
40	0.0027155221565254\\
45	0.0030535628757535\\
50	0.00339160488688012\\
55	0.00372964783270898\\
60	0.00406769137535546\\
65	0.0044057352575044\\
70	0.00474377928363676\\
75	0.00508182330085132\\
80	0.00541986718575898\\
85	0.00575791083585162\\
90	0.0060959541637962\\
95	0.00643399709440754\\
100	0.00677203958843614\\
};
\addlegendentry{MC$_8(0,0)$}

\addplot [color=red, line width=1.0pt, mark size=3.0pt, mark=diamond, mark options={solid, red}]
  table[row sep=crcr]{%
5	0.000196283169234979\\
10	0.000368634476511616\\
15	0.000539978785579487\\
20	0.000701724016912432\\
25	0.000858026187767255\\
30	0.00101297735382225\\
35	0.00116778286246687\\
40	0.00132265260637818\\
45	0.00147759572387453\\
50	0.0016325960304917\\
55	0.00178763933422802\\
60	0.00194271542576128\\
65	0.0020978170349327\\
70	0.00225293888124725\\
75	0.00240807704199341\\
80	0.00256322854278772\\
85	0.00271839108728229\\
90	0.00287356287425898\\
95	0.00302874247140488\\
100	0.00318392876872786\\
};
\addlegendentry{EC$_{10}$(0)}

\end{axis}

\end{tikzpicture}%
	\input{EC6clawerr.tex}
	\input{MC8clawerr.tex}
	\input{EC10clawerr.tex}
\end{center}
\caption{Solution and invariants error of different schemes.}
	\label{errorgrow}
\end{figure}

To show that the choice of parameters can affect solution accuracy significantly, the same problem is now solved with $c=5, T=5$, and $\Delta x=\Delta t=0.05$. The solution errors for $\mbox{MC}_6(\alpha_1)$, $\mbox{MC}_8(\alpha_2,\beta_2)$ and $\mbox{EC}_{10}(\alpha_3)$ are minimized by choosing $\alpha_1=8$, $(\alpha_2,\beta_2)=(-4,3.3)$ and $\alpha_3=-32$. 

Table~\ref{comp1solbbm} compares the schemes MC and EC introduced in Section~\ref{BBMmeth}, the LS and the PB schemes. 
All of the conservative schemes described in Section~\ref{BBMmeth} are able to preserve two discrete invariants (up to rounding errors). The solution error is at least comparable with the one given by the multisymplectic schemes but it is significantly smaller when choosing the optimal values of the free parameters. In particular, $\mbox{EC}_{10}(-32)$ is the most accurate scheme.

\begin{figure}[htbp]
\begin{center}
	\input{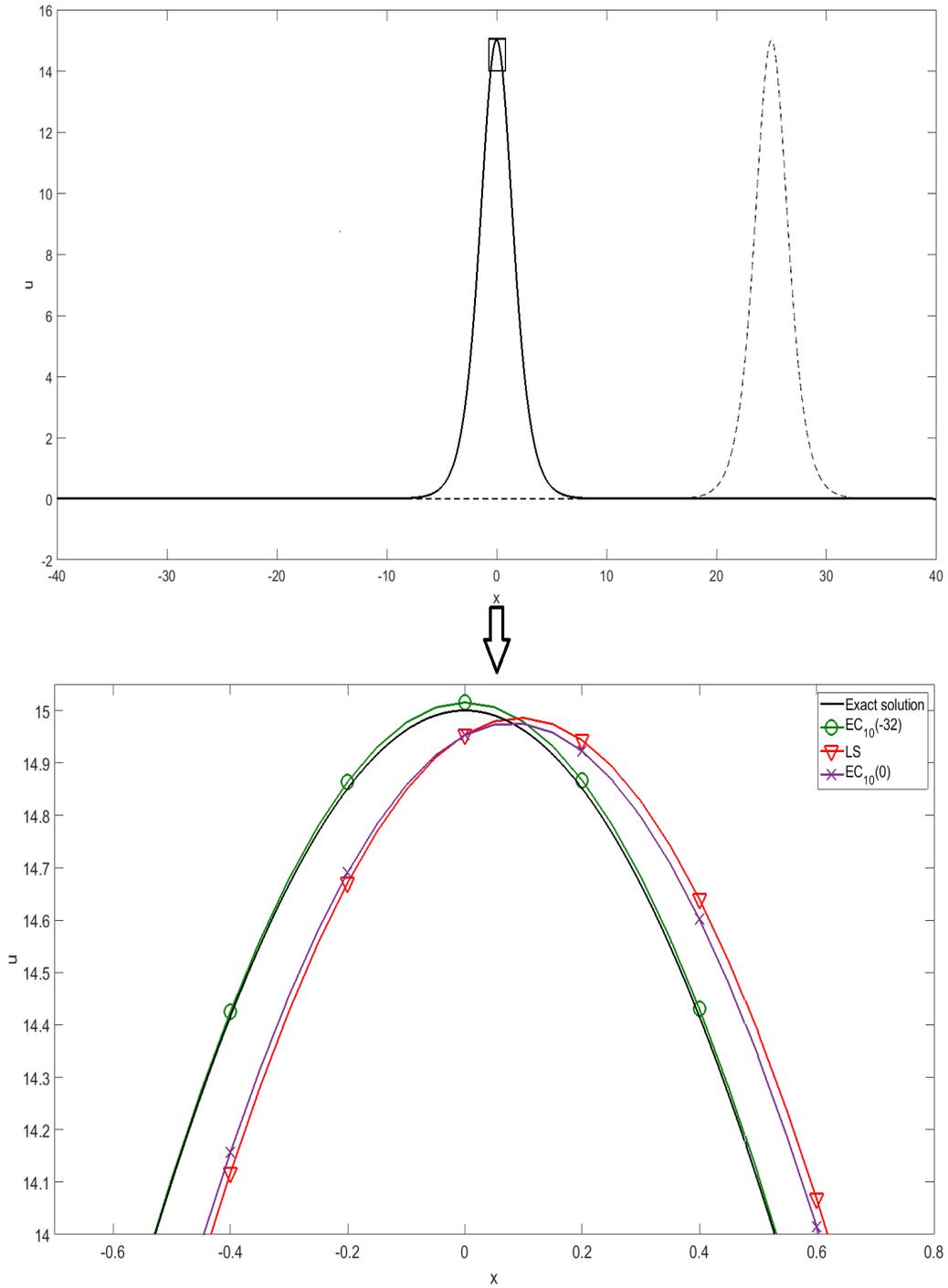}
\end{center}
\caption{Solitary wave problem for the BBM equation. Top: 
Initial condition (dashed line) and solution of $\mbox{EC}_{10}(-32)$ at time $T=5$, with $\Delta x=\Delta t=0.05$ (solid line). Bottom: Top of the single wave at time $T=5$; exact profile (solid line) and solutions of $\mbox{EC}_{10}(-32)$ (circles), LS (inverted triangles) and Koide \& Furihata scheme, $\mbox{EC}_{10}(0)$ (crosses). Markers are shown only at every fourth computed point.}
	\label{1solfigbbmnew}
\end{figure}

In Figure~\ref{1solfigbbmnew}, the upper plot shows the initial condition (dashed line) and the numerical solution given by method $\mbox{EC}_{10}(-32)$ at the final time $T=5$ (solid line). The lower plot shows the exact solution and the numerical approximations given by $\mbox{EC}_{10}(-32)$, LS and the Koide \& Furihata ($\mbox{EC}_{10}(0)$) schemes around the top of the wave. The solution of PB is not shown, as it is very close to the solution of LS. The scheme $\mbox{EC}_{10}(-32)$ gives the closest approximation to the exact solution. 

As a second benchmark problem, we study the interaction between two solitary waves. It is worth mentioning that the solitary wave (\ref{ex1solbbm}) is not a soliton; after interaction with another solitary wave, an oscillatory tail is generated and the solitary waves are not unscathed  (see, e.g., \cite{BonaPritchardScott,Olver}).

We consider equation (\ref{BBM}) with $\Omega=[-100,100]\times[0,15]$ and the initial condition
\begin{equation*}
u(x,0)=3c_1\sech^2\!\left(\tfrac{1}2(x-d_1)\right)+3c_2\sech^2\!\left(\tfrac{1}2(x-d_2)\right),
\end{equation*}
with
\begin{equation*}
c_1=6,\quad c_2=2,\quad d_1=40,\quad d_2=15.
\end{equation*}

The numerical solution of the Koide \& Furihata scheme ($\mbox{EC}_{10}(0)$) with $\Delta x=0.05$ and $\Delta t=0.003$ will be treated as the reference \textit{exact} solution\footnote{Similar results are obtained by using LS or PB on the same fine grid to compute a reference solution.}; we compare this with the solutions of the various schemes on a much coarser grid: $\Delta x=0.2$ and $\Delta t=0.015$.

The schemes $\mbox{MC}_{6}(\alpha_1)$, $\mbox{MC}_{8}(\alpha_2,\beta_2)$ and $\mbox{EC}_{10}(\alpha_3)$ each minimize their solution error when $\alpha_1=0.42$, $(\alpha_2,\beta_2)=(0.3,-0.04)$ and $\alpha_3=-0.4667$.

In Table~\ref{comp2solbbm} we compare the performance of the EC, MC and multisymplectic schemes, examining the errors in conservation laws, solution
and phase shift of the fastest wave at the final time $T=15$, defined as
$$\left.\text{Err}_\phi=(x_\text{max}-\tilde{x}_\text{max})\right\vert_{t=15};$$
here the locations of the peak of the fastest wave in the reference solution and the numerical solution are denoted by $x_\text{max}$ and $\tilde{x}_\text{max}$, respectively.

\begin{table}[t]
\caption{Solitary wave problem; $\Delta x=\Delta t=0.05$.} 
\label{comp1solbbm}
\small
\centerline{\begin{tabular}{|c|c|c|c|c|c|}
\hline
Method &  $\text{Err}_1$ & $\text{Err}_2$ & $\text{Err}_3$  &  Error in the solution\\ 
\hline
$\mbox{EC}_6$	& 1.49e-13   & 0.0035 & 8.64e-12 &  0.0375\\	
\hline   
$\mbox{MC}_6(0)$	& 2.06e-13 & 1.19e-12 & 0.0199 & 0.0445\\  
\hline
$\mbox{MC}_6(8)$	& 2.42e-13  & 1.02e-12 & 0.0428 & 0.0062 \\	
\hline
$\mbox{EC}_{8}$	& 2.27e-13  & 0.0042 & 7.73e-12 &  0.0375\\
\hline
$\mbox{MC}_{8}(0,0)$	& 2.91e-13  & 1.76e-12 & 0.0162 &  0.0443\\	
\hline
$\mbox{MC}_{8}(-4,3.3)$	& 2.13e-13  & 1.65e-12 & 0.1857 & 0.0043 \\	
\hline
Koide \& Furihata; $\mbox{EC}_{10}(0)$	& 1.71e-13  & 0.0034 & 7.28e-12 &  0.0356\\	
\hline
$\mbox{EC}_{10}(-32)$	& 2.34e-13  & 0.0030& 6.82e-12 & 0.0037 \\
\hline
LS & 1.99e-13  & 2.63e-04 & 0.0100 &  0.0415\\
\hline
PB & 1.49e-13  & 0.1379 & 1.3841 &  0.0434\\
\hline
\end{tabular}}
\end{table}

\begin{table}[htb]
\caption{Two solitary wave problem for the BBM equation; $\Delta x=0.2$, $\Delta t=0.015$.}
\label{comp2solbbm}
\small
\centerline{\begin{tabular}{|c|c|c|c|c|c|}
\hline
Method &  $\text{Err}_1$ & $\text{Err}_2$ & $\text{Err}_3$  &  Solution error & $\text{Err}_\phi$ \\
\hline
$\mbox{EC}_6$	& 2.27e-13   & 0.2270 & 1.55e-11 &  0.1251 & -0.3\\
\hline   
$\mbox{MC}_6(0)$	& 2.70e-13 & 2.96e-12 & 7.8731 & 0.1112& -0.3\\
\hline
$\mbox{MC}_6(0.42)$	& 2.56e-13  & 1.82e-12 & 7.7930 & 0.0038& 0.1\\
\hline
$\mbox{EC}_{8}$	& 3.41e-13  & 0.6306 & 1.91e-11 &  0.1251& -0.3\\
\hline
$\mbox{MC}_{8}(0,0)$	& 2.56e-13  & 2.16e-12 & 4.1856 &  0.1081& -0.3\\	
\hline
$\mbox{MC}_{8}(0.3,-0.04)$	& 3.55e-13  & 3.52e-12 & 7.5081 & 0.0097 & 0.1\\	
\hline
Koide \& Furihata; $\mbox{EC}_{10}(0)$	& 2.70e-13  & 0.2467 & 1.55e-11 &  0.0256 & -0.1\\
\hline
$\mbox{EC}_{10}(-0.4667)$	& 2.98e-13  & 0.2463 & 2.00e-11 & 0.0125 & 0.1\\
\hline
LS	&  2.98e-13  & 0.5514 & 9.8587 & 0.0398  & 0.1 \\
\hline
PB	& 1.28e-13  & 0.0983 & 0.6815 &  0.0630 & -0.1\\ 
\hline
\end{tabular}}
\end{table}

Table~\ref{comp2solbbm} shows that each $\mbox{MC}$ and $\mbox{EC}$ scheme preserves the discretization of two invariants in (\ref{glinvbbm}) to machine accuracy. The scheme $\mbox{MC}_6(0.42)$ gives the smallest solution error.

\begin{figure}[htbp]
{\includegraphics[width=16.5cm,height=9cm]{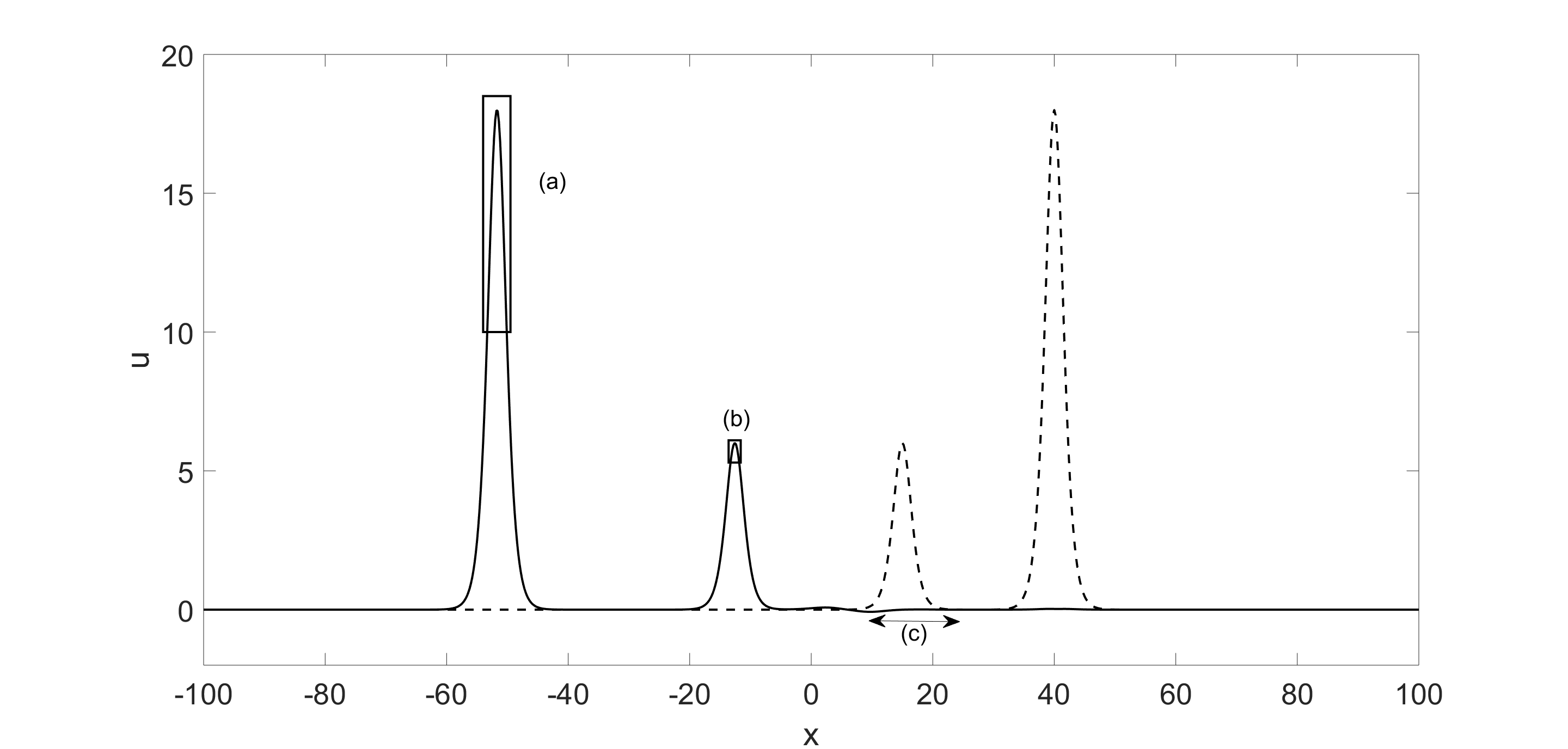}}
\begin{center}
	\definecolor{mycolor1}{rgb}{0.00000,0.49804,0.00000}%
\definecolor{mycolor2}{rgb}{0.49412,0.18431,0.55686}%
\begin{tikzpicture}

\begin{axis}[%
width=5in,
height=3in,
at={(1in,7in)},
scale only axis,
xmin=-54,
xmax=-49.5,
xlabel style={font=\color{white!15!black}},
xlabel={$x$},
ymin=10,
ymax=18.5,
ylabel style={font=\color{white!15!black}},
ylabel={$u$},
axis background/.style={fill=white},
legend style={legend cell align=left, align=left, draw=white!15!black}
]
\addplot [color=black, line width=2.0pt]
  table[row sep=crcr]{%
-55	2.50839814774226\\
-54.95	2.62722303530553\\
-54.9	2.75117568154351\\
-54.85	2.88042704732892\\
-54.8	3.01514859004679\\
-54.75	3.15551162138627\\
-54.7	3.3016865898583\\
-54.65	3.45384228330484\\
-54.6	3.61214494673605\\
-54.55	3.7767573109725\\
-54.5	3.94783752779398\\
-54.45	4.12553800761486\\
-54.4	4.31000415613224\\
-54.35	4.50137300694064\\
-54.3	4.69977174778936\\
-54.25	4.90531613898993\\
-54.2	5.11810882347526\\
-54.15	5.33823752918204\\
-54.1	5.56577316578683\\
-54.05	5.8007678193831\\
-54	6.04325265045171\\
-53.95	6.29323570245504\\
-53.9	6.5506996305778\\
-53.85	6.81559936254225\\
-53.8	7.08785970603477\\
-53.75	7.36737292007885\\
-53.7	7.65399627065459\\
-53.65	7.94754959396691\\
-53.6	8.24781289396279\\
-53.55	8.55452400394244\\
-53.5	8.86737634533959\\
-53.45	9.18601681988962\\
-53.4	9.510043874378\\
-53.35	9.83900577987188\\
-53.3	10.1723991696793\\
-53.25	10.5096678821427\\
-53.2	10.8502021556329\\
-53.15	11.1933382236514\\
-53.1	11.538358357633\\
-53.05	11.8844914037692\\
-53	12.2309138578087\\
-52.95	12.5767515182568\\
-52.9	12.9210817535944\\
-52.85	13.2629364130286\\
-52.8	13.6013054028373\\
-52.75	13.9351409416004\\
-52.7	14.263362497576\\
-52.65	14.5848624002818\\
-52.6	14.8985121061415\\
-52.55	15.2031690850493\\
-52.5	15.497684281153\\
-52.45	15.7809100873658\\
-52.4	16.0517087594211\\
-52.35	16.3089611820764\\
-52.3	16.5515758877408\\
-52.25	16.7784982167717\\
-52.2	16.9887194993451\\
-52.15	17.1812861315642\\
-52.1	17.3553084136482\\
-52.05	17.509969015954\\
-52	17.6445309394336\\
-51.95	17.7583448410582\\
-51.9	17.8508556017863\\
-51.85	17.9216080247588\\
-51.8	17.9702515643892\\
-51.75	17.9965440026267\\
-51.7	18.0003540065172\\
-51.65	17.9816625208391\\
-51.6	17.9405629705057\\
-51.55	17.8772602690486\\
-51.5	17.7920686512204\\
-51.45	17.6854083689743\\
-51.4	17.5578013102203\\
-51.35	17.4098656182695\\
-51.3	17.2423094062914\\
-51.25	17.0559236750166\\
-51.2	16.8515745530243\\
-51.15	16.6301949870454\\
-51.1	16.392776014698\\
-51.05	16.1403577539528\\
-51	15.8740202425112\\
-50.95	15.5948742563527\\
-50.9	15.3040522302478\\
-50.85	15.002699394365\\
-50.8	14.6919652305943\\
-50.75	14.3729953402761\\
-50.7	14.046923802064\\
-50.65	13.7148660850776\\
-50.6	13.3779125687036\\
-50.55	13.0371227067292\\
-50.5	12.6935198602759\\
-50.45	12.3480868114974\\
-50.4	12.0017619584517\\
-50.35	11.6554361811134\\
-50.3	11.309950359285\\
-50.25	10.9660935152694\\
-50.2	10.6246015476001\\
-50.15	10.2861565168881\\
-50.1	9.95138644088942\\
-50.05	9.62086555314771\\
-50	9.29511497793385\\
-49.95	8.97460377358287\\
-49.9	8.65975029658998\\
-49.85	8.3509238398602\\
-49.8	8.04844650017647\\
-49.75	7.75259523214464\\
-49.7	7.46360404847084\\
-49.65	7.18166632932199\\
-49.6	6.90693720661721\\
-49.55	6.63953599230646\\
-49.5	6.37954862294013\\
-49.45	6.12703009605376\\
-49.4	5.88200687703086\\
-49.35	5.64447925812278\\
-49.3	5.41442365416146\\
-49.25	5.19179482217692\\
-49.2	4.97652799460535\\
-49.15	4.76854091803986\\
-49.1	4.56773579152562\\
-49.05	4.37400110023862\\
-49	4.18721334201292\\
};
\addlegendentry{Exact solution}

\addplot [color=mycolor1, line width=1.0pt, only marks, mark size=4.0pt, mark=o, mark options={solid, mycolor1}, mark repeat=2]
  table[row sep=crcr]{%
-55	2.47513009050238\\
-54.8	2.97936641341105\\
-54.6	3.57422532531686\\
-54.4	4.27051536419028\\
-54.2	5.07782847224275\\
-54	6.00316854402039\\
-53.8	7.04914585607551\\
-53.6	8.21177401637858\\
-53.4	9.47802574769907\\
-53.2	10.8234723170212\\
-53	12.2105261361868\\
-52.8	13.5879697995427\\
-52.6	14.8924943775323\\
-52.4	16.0527785576748\\
-52.2	16.9961529200576\\
-52	17.6571590326732\\
-51.8	17.9865394374653\\
-51.6	17.9587088040247\\
-51.4	17.5758556593824\\
-51.2	16.8675858385498\\
-51	15.8862041783343\\
-50.8	14.6988779406036\\
-50.6	13.3785926912242\\
-50.4	11.9958077894013\\
-50.2	10.6121660680449\\
-50	9.27682511580124\\
-49.8	8.0252716053439\\
-49.6	6.88003895681594\\
-49.4	5.85259910915477\\
-49.2	4.94576496708918\\
-49	4.15611389354822\\
};
\addlegendentry{MC$_6(0.42)$}
\addplot [color=red, line width=1.0pt, mark size=4.0pt, mark=triangle, mark options={solid, rotate=180, red}, mark repeat=2, mark phase=2]
  table[row sep=crcr]{%
-55	2.68669436455476\\
-54.8	3.22616665761405\\
-54.6	3.86036416757167\\
-54.4	4.59980941070187\\
-54.2	5.4533984893037\\
-54	6.42687721617246\\
-53.8	7.52086880421963\\
-53.6	8.72851426811362\\
-53.4	10.0329257951742\\
-53.2	11.4048447415701\\
-53	12.8011095614414\\
-52.8	14.1647029184617\\
-52.6	15.4271473714566\\
-52.4	16.5137309770161\\
-52.2	17.3514071399734\\
-52	17.8783240350509\\
-51.8	18.0530989475415\\
-51.6	17.861582438054\\
-51.4	17.3192682469241\\
-51.2	16.4686569593791\\
-51	15.3723295260905\\
-50.8	14.1036118697508\\
-50.6	12.7370844139195\\
-50.4	11.3407895978801\\
-50.2	9.97113921284052\\
-50	8.67064087162974\\
-49.8	7.46794136510169\\
-49.6	6.37941205111719\\
-49.4	5.41151196779203\\
-49.2	4.56333283062569\\
-49	3.82894392432891\\
};
\addlegendentry{LS}
\addplot [color=blue, line width=1.0pt, mark size=4.0pt, mark=triangle, mark options={solid, blue}, mark repeat=2, mark phase=2]
  table[row sep=crcr]{%
-55	2.17461365627983\\
-54.8	2.62122831192529\\
-54.6	3.15025535246054\\
-54.4	3.7726534353639\\
-54.2	4.49888818828895\\
-54	5.33786268683182\\
-53.8	6.29543953140091\\
-53.6	7.37253640066613\\
-53.4	8.56285948101982\\
-53.2	9.85047102147825\\
-53	11.2075637581348\\
-52.8	12.5930055567368\\
-52.6	13.9523558821416\\
-52.4	15.2200428018256\\
-52.2	16.324124461269\\
-52	17.1935009448693\\
-51.8	17.7666793758619\\
-51.6	18.0004725399129\\
-51.4	17.8766626726492\\
-51.2	17.4049407797954\\
-51	16.6213338218244\\
-50.8	15.5825442693024\\
-50.6	14.3576701355575\\
-50.4	13.0192607454514\\
-50.2	11.6354910127991\\
-50	10.2645943004105\\
-49.8	8.95191065003306\\
-49.6	7.72926878221983\\
-49.4	6.61606439384578\\
-49.2	5.62131588332881\\
-49	4.7460816684868\\
};
\addlegendentry{PB}
\addplot [color=mycolor2, line width=1.0pt, mark size=4.0pt, mark=x, mark options={solid, mycolor2}, mark repeat=2]
  table[row sep=crcr]{%
-55	2.37579796736869\\
-54.8	2.85586704029281\\
-54.6	3.42235951050201\\
-54.4	4.08605752865602\\
-54.2	4.85690633497755\\
-54	5.74281193839424\\
-53.8	6.74801402929485\\
-53.6	7.87103336911038\\
-53.4	9.10229040290233\\
-53.2	10.421642389607\\
-53	11.7962809661716\\
-52.8	13.1796307068657\\
-52.6	14.5120102140994\\
-52.4	15.7237447347589\\
-52.2	16.741045167441\\
-52	17.4942712321786\\
-51.8	17.9273245354429\\
-51.6	18.0062016609392\\
-51.4	17.7245713952866\\
-51.2	17.1048444325781\\
-51	16.1944306783994\\
-50.8	15.0582355620037\\
-50.6	13.7693639669154\\
-50.4	12.4001612028463\\
-50.2	11.0151997850871\\
-50	9.66696795598781\\
-49.8	8.3942127297247\\
-49.6	7.22236984537837\\
-49.4	6.16531959734416\\
-49.2	5.22775951949179\\
-49	4.40766229569852\\
};
\addlegendentry{EC$_{10}(0)$}
\end{axis}
\end{tikzpicture}%
\end{center}
\caption{Two-wave problem for the BBM equation. Top: Initial condition (dashed line) and solution of $\mbox{MC}_6(0.42)$ at time $T=15$ (solid line). Bottom: Magnification around zone (a); Reference solution (bold line) and solutions of $\mbox{MC}_6(0.42)$ (circles), LS (inverted triangles), PB (triangles) and $\mbox{EC}_{10}(0)$ (crosses); markers at every second point.}\label{ft2BBM}
\end{figure}

\begin{figure}[htbp]
\begin{center}
	\definecolor{mycolor1}{rgb}{0.00000,0.49804,0.00000}%
\definecolor{mycolor2}{rgb}{0.49412,0.18431,0.55686}%
\begin{tikzpicture}

\begin{axis}[%
width=5in,
height=3in,
at={(1in,7in)},
scale only axis,
xmin=-13.6,
xmax=-11.6,
xlabel style={font=\color{white!15!black}},
xlabel={$x$},
ymin=5.3,
ymax=6.1,
ylabel style={font=\color{white!15!black}},
ylabel={$u$},
axis background/.style={fill=white},
legend style={legend cell align=left, align=left, draw=white!15!black}
]
\addplot [color=black, line width=2.0pt]
  table[row sep=crcr]{%
-15	1.76463547152339\\
-14.95	1.83998299299281\\
-14.9	1.91781323916199\\
-14.85	1.99813797954417\\
-14.8	2.08096142185352\\
-14.75	2.16627945992523\\
-14.7	2.25407889817333\\
-14.65	2.34433665725914\\
-14.6	2.4370189665649\\
-14.55	2.53208055004465\\
-14.5	2.62946381305044\\
-14.45	2.72909803879404\\
-14.4	2.8308986041823\\
-14.35	2.93476622584583\\
-14.3	3.04058624823318\\
-14.25	3.14822798665068\\
-14.2	3.25754413904658\\
-14.15	3.36837028114864\\
-14.1	3.4805244602165\\
-14.05	3.59380690313227\\
-14	3.70799985478204\\
-13.95	3.82286756263399\\
-13.9	3.93815642305899\\
-13.85	4.05359530422304\\
-13.8	4.16889605928131\\
-13.75	4.28375424208431\\
-13.7	4.39785003565558\\
-13.65	4.51084940130316\\
-13.6	4.62240545338608\\
-13.55	4.73216006149183\\
-13.5	4.83974567811884\\
-13.45	4.94478738595017\\
-13.4	5.04690515451459\\
-13.35	5.14571629154181\\
-13.3	5.2408380697282\\
-13.25	5.33189050504902\\
-13.2	5.41849925830922\\
-13.15	5.50029862744931\\
-13.1	5.57693459435153\\
-13.05	5.64806788666265\\
-13	5.71337701259182\\
-12.95	5.77256122487423\\
-12.9	5.82534336921489\\
-12.85	5.87147257261842\\
-12.8	5.91072672812426\\
-12.75	5.94291473461769\\
-12.7	5.96787845356455\\
-12.65	5.98549434867076\\
-12.6	5.99567477951123\\
-12.55	5.99836892599319\\
-12.5	5.99356332696599\\
-12.45	5.98128202319637\\
-12.4	5.96158630210634\\
-12.35	5.93457404891955\\
-12.3	5.90037871598374\\
-12.25	5.85916792883194\\
-12.2	5.81114175383216\\
-12.15	5.75653065788862\\
-12.1	5.69559319546173\\
-12.05	5.62861346205712\\
-12	5.55589835622641\\
-11.95	5.47777469397934\\
-11.9	5.39458622032879\\
-11.85	5.30669056250231\\
-11.8	5.21445616821564\\
-11.75	5.11825927040324\\
-11.7	5.0184809170392\\
-11.65	4.91550410128758\\
-11.6	4.80971102332139\\
-11.55	4.70148051088407\\
-11.5	4.59118562117194\\
-11.45	4.47919144202135\\
-11.4	4.36585310581024\\
-11.35	4.25151402504309\\
-11.3	4.13650435436878\\
-11.25	4.02113967986647\\
-11.2	3.90571993287929\\
-11.15	3.79052852252824\\
-11.1	3.67583167832215\\
-11.05	3.56187799200889\\
-11	3.44889814598425\\
-10.95	3.33710481417832\\
-10.9	3.22669272034904\\
-10.85	3.11783883810274\\
-10.8	3.01070271669106\\
-10.75	2.90542691667177\\
-10.7	2.80213753981871\\
-10.65	2.70094483818516\\
-10.6	2.60194388792828\\
-10.55	2.50521531434395\\
-10.5	2.4108260555107\\
-10.45	2.31883015296341\\
-10.4	2.22926955888178\\
-10.35	2.14217495036013\\
-10.3	2.05756654240184\\
-10.25	1.97545489233389\\
-10.2	1.89584168935049\\
-10.15	1.81872052385852\\
-10.1	1.7440776322001\\
-10.05	1.67189261316542\\
-10	1.6021391134781\\
-9.95	1.5347854801307\\
-9.90000000000001	1.46979537807457\\
-9.84999999999999	1.407128372323\\
-9.8	1.34674047401291\\
-9.75	1.28858465039481\\
-9.7	1.23261129908008\\
-9.65000000000001	1.17876868718136\\
-9.59999999999999	1.12700335623378\\
-9.55	1.07726049399116\\
-9.5	1.02948427435285\\
-9.45	0.983618166802444\\
-9.40000000000001	0.939605216830084\\
-9.34999999999999	0.897388298870647\\
-9.3	0.856910343325811\\
-9.25	0.818114539252103\\
-9.2	0.78094451429007\\
-9.15000000000001	0.745344493390724\\
-9.09999999999999	0.71125943786037\\
-9.05	0.678635166201268\\
-9	0.647418458173541\\
-8.95	0.617557143445253\\
-8.90000000000001	0.589000176133329\\
-8.84999999999999	0.561697696474753\\
-8.8	0.535601080795736\\
-8.75	0.510662980880564\\
-8.7	0.486837353772888\\
-8.65000000000001	0.464079482971928\\
-8.59999999999999	0.442345991923329\\
-8.55	0.421594850636618\\
-8.5	0.401785376200662\\
-8.45	0.382878227908226\\
-8.40000000000001	0.364835397643831\\
-8.34999999999999	0.347620196134557\\
-8.3	0.331197235613177\\
-8.25	0.315532409393437\\
-8.2	0.300592868813542\\
-8.15000000000001	0.286346997960066\\
-8.09999999999999	0.272764386546706\\
-8.05	0.259815801284286\\
-8	0.247473156044848\\
};
\addlegendentry{Exact solution}

\addplot [color=mycolor1, line width=1.0pt, only marks, mark size=4.0pt, mark=o, mark options={solid, mycolor1}]
  table[row sep=crcr]{%
-15	1.76859145770486\\
-14.8	2.08746387794617\\
-14.6	2.44645763778838\\
-14.4	2.84354231209646\\
-14.2	3.27345840316124\\
-14	3.72696170619365\\
-13.8	4.19032372784445\\
-13.6	4.64532551344635\\
-13.4	5.06997894502258\\
-13.2	5.44012105730589\\
-13	5.7318387558835\\
-12.8	5.92442300614225\\
-12.6	6.00330611129355\\
-12.4	5.96231627193247\\
-12.2	5.80467658624152\\
-12	5.5424801706901\\
-11.8	5.19478342279199\\
-11.6	4.78481330020735\\
-11.4	4.33695003238028\\
-11.2	3.87408854909312\\
-11	3.41576416660448\\
-10.8	2.97716257769546\\
-10.6	2.56891793499311\\
-10.4	2.19748254016112\\
-10.2	1.86582454714147\\
-10	1.57424546439794\\
-9.8	1.3211716477877\\
-9.59999999999999	1.10383697733836\\
-9.40000000000001	0.91882367150581\\
-9.2	0.762460959307136\\
-9	0.631099125742467\\
-8.8	0.521283449312965\\
-8.59999999999999	0.429852931448111\\
-8.40000000000001	0.353985708841735\\
-8.2	0.291208774660944\\
-8	0.239385333788531\\
};
\addlegendentry{MC$_6(0.42)$}
\addplot [color=red, line width=1.0pt, mark size=4.0pt, mark=triangle, mark options={solid, rotate=180, red}]
  table[row sep=crcr]{%
-15	1.87578362898547\\
-14.8	2.20792097526141\\
-14.6	2.5800282536254\\
-14.4	2.98925433501363\\
-14.2	3.42924415149613\\
-14	3.88938941803107\\
-13.8	4.35439473772199\\
-13.6	4.80442193624522\\
-13.4	5.21606062379307\\
-13.2	5.56425518288346\\
-13	5.82508466594668\\
-12.8	5.97898751318965\\
-12.6	6.01376000569002\\
-12.4	5.92657737515552\\
-12.2	5.72447731073143\\
-12	5.42316616933122\\
-11.8	5.04448894800104\\
-11.6	4.61323982189456\\
-11.4	4.15406016221132\\
-11.2	3.68899575602612\\
-11	3.2359867862234\\
-10.8	2.8082800854743\\
-10.6	2.41456805645694\\
-10.4	2.05958929844665\\
-10.2	1.74494370838353\\
-10	1.46993695558463\\
-9.8	1.23234054183186\\
-9.59999999999999	1.02901383006452\\
-9.40000000000001	0.856376243443212\\
-9.2	0.710742142069616\\
-9	0.588541935596146\\
-8.8	0.486455528248208\\
-8.59999999999999	0.401482087534566\\
-8.40000000000001	0.330966032443819\\
-8.2	0.272594633258751\\
-8	0.224378512269897\\
};
\addlegendentry{LS}
\addplot [color=blue, line width=1.0pt, mark size=4.0pt, mark=triangle, mark options={solid, blue}]
  table[row sep=crcr]{%
-15	1.68688158172296\\
-14.8	1.9935049405236\\
-14.6	2.3401030231818\\
-14.4	2.72542223570046\\
-14.2	3.14524694636201\\
-14	3.59161697612675\\
-13.8	4.05223631215692\\
-13.6	4.51029676639295\\
-13.4	4.94495604077357\\
-13.2	5.33265073232776\\
-13	5.64926830175246\\
-12.8	5.87296160845327\\
-12.6	5.98713169263669\\
-12.4	5.98293653932747\\
-12.2	5.86070633937746\\
-12	5.62988968811913\\
-11.8	5.30754531492931\\
-11.6	4.91577826215037\\
-11.4	4.47874903707585\\
-11.2	4.01989151153434\\
-11	3.5597973451156\\
-10.8	3.11496479471836\\
-10.6	2.69737371420794\\
-10.4	2.31469900630075\\
-10.2	1.97092213050105\\
-10	1.6671199398756\\
-9.8	1.40226666018473\\
-9.59999999999999	1.17394870587682\\
-9.40000000000001	0.978945892198569\\
-9.2	0.81367011949038\\
-9	0.674474422183781\\
-8.8	0.557855118197199\\
-8.59999999999999	0.46057189937841\\
-8.40000000000001	0.379708573414917\\
-8.2	0.312693228760382\\
-8	0.257292303963527\\
};
\addlegendentry{PB}
\addplot [color=mycolor2, line width=1.0pt, only marks, mark size=4.0pt, mark=x, mark options={solid, mycolor2}]
  table[row sep=crcr]{%
-15	1.78226983410182\\
-14.8	2.10020424093472\\
-14.6	2.4579314674957\\
-14.4	2.85352491639665\\
-14.2	3.28188157625445\\
-14	3.73394401972057\\
-13.8	4.19615680026882\\
-13.6	4.6504003722398\\
-13.4	5.07465481739007\\
-13.2	5.44456749258172\\
-13	5.73590940823867\\
-12.8	5.9276270581557\\
-12.6	6.00491648626543\\
-12.4	5.9615964793221\\
-12.2	5.80114802091085\\
-12	5.53612547791336\\
-11.8	5.18610758333614\\
-11.6	4.77474735657141\\
-11.4	4.32664570114473\\
-11.2	3.86468346901842\\
-11	3.4081896598795\\
-10.8	2.97203324311475\\
-10.6	2.56650668588892\\
-10.4	2.19776059915311\\
-10.2	1.868537223348\\
-10	1.57899746014614\\
-9.8	1.32750402647858\\
-9.59999999999999	1.11128683989577\\
-9.40000000000001	0.926964440622354\\
-9.2	0.770925301122029\\
-9	0.639588364522149\\
-8.8	0.529567622471699\\
-8.59999999999999	0.437765105467751\\
-8.40000000000001	0.361413341324036\\
-8.2	0.298084063239829\\
-8	0.245675779917967\\
};
\addlegendentry{EC$_{10}(0)$}
\end{axis}
\end{tikzpicture}%
	\input{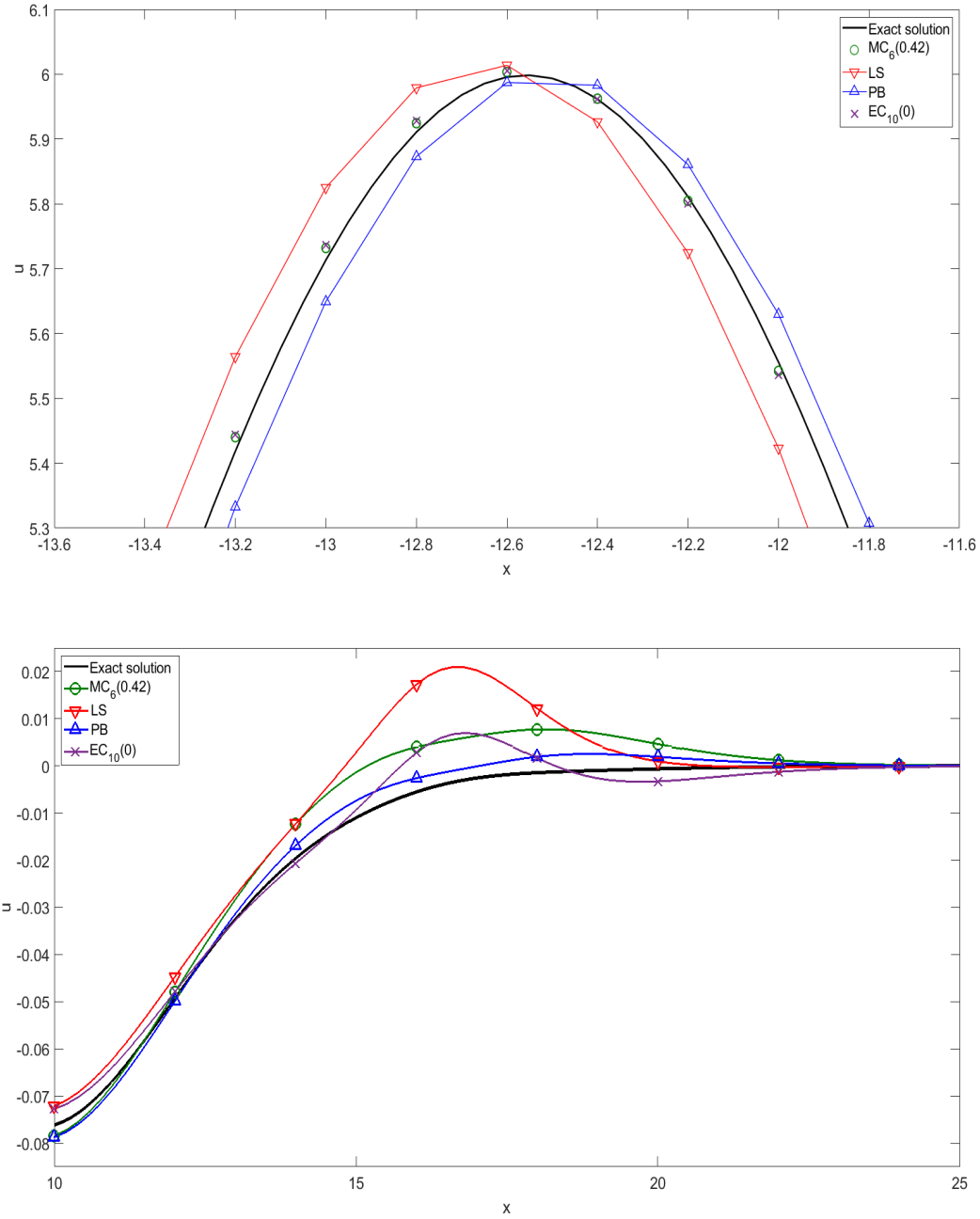}
\end{center}
\caption{Two-wave problem for the BBM equation. Reference solution (bold line) and solutions of $\mbox{MC}_6(0.42)$ (circles), LS (inverted triangles), PB (triangles) and $\mbox{EC}_{10}(0)$ (crosses) at time $T=15$. Top: magnification of zone (b) in Figure~\ref{ft2BBM}. Bottom: magnification of zone (c) in Figure~\ref{ft2BBM}; markers are at every tenth point.}\label{solitBBM}
\end{figure}

The upper part of Figure~\ref{ft2BBM} shows the initial condition (dashed line) and the numerical solution given by method $\mbox{MC}_6(0.42)$ at time $T=15$ (solid line). The lower part of Figure~\ref{ft2BBM} and the upper part of Figure~\ref{solitBBM} show the reference solution and the numerical solutions given by MC$_6$(0.42), LS, PB and the scheme of Koide \& Furihata (EC$_{10}$(0)) around the top of the two waves. The zone around the top of the faster wave is where the most obvious differences between schemes occur. The two waves are approximated best overall by MC$_6$(0.42), but PB gives a slightly better approximation of the oscillatory tail generated after the interaction of the two waves (see the lower part of Figure~\ref{solitBBM}).

\section{Nonlinear Schr\"odinger equation}\label{NLSsec}

In this section we develop schemes for the nonlinear Schr\"odinger (NLS) equation with complex variable $\psi$:
\begin{equation}\label{NLSc}
i\psi_t+\psi_{xx}+|\psi|^2\psi=0,\quad(x,t)\in\Omega\equiv[a,b]\times[0,\infty).
\end{equation}
Setting $\psi=u+iv$, with $u,v\in\mathbb{R}$, (\ref{NLSc}) reduces to the following system of equations:
\begin{equation}\label{NLS}
\mca=\left(A[u,v],\,A[-v,u]\right)=\mathbf{0}, \qquad
A[a,b]\equiv a_t+b_{xx}+(a^2+b^2)b,\qquad (x,t)\in\Omega.
\end{equation}
The NLS equation has infinitely many independent conservation laws \cite{FT87}. The first three are:
\begin{align}\label{NLSCL1}
&\mca\mcq_1= D_x(F_1)+D_t(G_1) = D_x\left(2uv_x-2u_xv\right)+D_t(u^2+v^2)=0,\\\label{NLSCL2}
&\mca\mcq_2= D_x(F_2)\!+\!D_t(G_2) = D_x\!\left(u_x^2+v_x^2+u_tv-uv_t+\tfrac{1}2(u^2+v^2)^2\right)+D_t(uv_x\!-\!u_xv)=0,\\\label{NLSCL3}
&\mca\mcq_3= D_x(F_3)+D_t(G_3) = D_x\left(-2u_xu_t-2v_tv_x\right)+D_t\left(u_x^2+v_x^2-\tfrac{1}2(u^2+v^2)^2\right)=0;
\end{align}
their characteristics are, respectively,
$$\mcq_1=(2u,-2v)^T,\quad \mcq_2=(2v_x,2u_x)^T,\quad \mcq_3=(-2v_t,-2u_t)^T.$$
Integrating these conservation laws in space, subject to conditions that give no boundary contributions, one obtains the global invariants
$$\int\left( u^2+v^2\right)\mathrm{d}x,\quad\int\left( uv_x-u_xv\right)\mathrm{d}x,\quad\int\left( u_x^2+v_x^2-\tfrac{1}2{(u^2+v^2)^2}\right)\mathrm{d}x,$$
which represent charge, momentum and energy, respectively \cite{FT87}. The NLS system (\ref{NLS}) can be written in Hamiltonian form,
$$\left(\begin{array}{c}\mathrm{d}u/\mathrm{d}t\\\mathrm{d}v/\mathrm{d}t\end{array}\right)=\left(\begin{array}{cc}
0&1\\
-1&0
\end{array}\right) \left(\begin{array}{c} \delta\mathcal{H}/\delta u \\\delta\mathcal{H}/\delta v \end{array}\right),$$
with the Hamiltonian functional
\begin{equation}\label{Hamiltonian}
\mathcal{H}=\tfrac{1}{2}\int G_3\,\mathrm{d}x.
\end{equation}

\subsection{Compact conservative methods for the NLS equation}
In this section we develop compact methods that preserve multiple local conservation laws of the NLS equation. We choose $A=-1$, $B=1$ in Figure~\ref{stencil} and seek second-order approximations of characteristics, densities and fluxes at $(0,1/2)$, $(0,0)$ and $(-1/2,1/2)$, respectively.

The schemes obtained using the strategy described in Section~\ref{method} depend on some free parameters that are $\mathcal{O}(\Delta x^2,\Delta t^2)$. The free parameters can be chosen to minimize the solution error for a given  problem, but no choice gives a higher-order scheme. 

\subsubsection*{3-parameter family of local energy-conserving methods}
Charge- and energy-conserving schemes are found by specifying the approximations of the higher-order derivatives and using compact discretizations of the characteristics and the linear factors of the nonlinear terms in (\ref{NLS}). Hence, we consider discretizations of the form
\[ 
\widetilde{\mca}=(\widetilde{A}_1[u,v],\,\widetilde{A}_2[u,v])=\mathbf{0},
\]
where
\[
\widetilde{A}_1[u,v]=\widetilde{u_t}+D_m^2\mu_nv_{-1,0}+(\widetilde{u^2}+\widetilde{v^2})\mu_nv_{0,0},\quad \widetilde{A}_2[u,v]=-\widetilde{v_t}+D_m^2\mu_nu_{-1,0}+(\widehat{u^2}+\widehat{v^2})\mu_nu_{0,0},
\]
$\widetilde{u_t}$ and $\widetilde{v_t}$ are of the form (\ref{linapp}) and $\widetilde{u^2}$, $\widetilde{v^2}$, $\widehat{u^2}$ and $\widehat{v^2}$ are of the form (\ref{quadapp}). The discrete characteristics are
$$\widetilde{\mcq}_1=(2\mu_nu_{0,0},-2\mu_nv_{0,0})^T,\qquad \widetilde{\mcq}_3=(-2D_nv_{0,0},-2D_nu_{0,0})^T.$$
Solving the conditions
$$\mathsf{E}(\widetilde{\mca}\widetilde{\mcq}_1)=\mathbf{0},\qquad\mathsf{E}(\widetilde{\mca}\widetilde{\mcq}_3)=\mathbf{0}$$
yields a 3-parameter family of methods, defined by
$$\widetilde{\mca}(\lambda,\eta,\nu)=(\widetilde{A}[u,v],\,\widetilde{A}[-v,u])^T,$$
where
\begin{align*}
\widetilde{A}[a,b]=&\,D_n(a_{0,0}+\lambda D_m^2a_{-1,0})+D_m^2\mu_nb_{-1,0}\\
&+\left[\mu_n(a_{0,0}^2+b_{0,0}^2)+\eta D_m^2\mu_n(a_{-1,0}^2+b_{-1,0}^2)+\nu D_mD_n\mu_m(a_{-1,0}^2+b_{-1,0}^2)\right]\mu_nb_{0,0},
\end{align*}
and $\lambda$, $\eta$ and $\nu$ are $\mathcal{O}(\Delta x^2,\Delta t^2)$. 
The scheme $\widetilde{\mca}(\lambda,\eta,\nu)$ has the following discrete local conservation laws for charge and energy:
\begin{equation*}
 \widetilde{\mca}\widetilde{\mcq}_1=D_m\widetilde{F}_1+D_n\widetilde{G}_1=0,\qquad\widetilde{\mca}\widetilde{\mcq}_3=D_m\widetilde{F}_3+D_n\widetilde{G}_3=0,
\end{equation*}
where
\begin{align}\nonumber
\widetilde{F}_1=&\,2(\mu_m\mu_nu_{-1,0})(D_m\mu_nv_{-1,0})-2(D_m\mu_nu_{-1,0})(\mu_m\mu_nv_{-1,0})+\lambda(\Theta[u,u]+\Theta[v,v]),\\\label{dfec}
\widetilde{G}_1=&\,u_{0,0}^2+v_{0,0}^2+\lambda(u_{0,0}D_m^2u_{-1,0}+v_{0,0}D_m^2v_{-1,0}),\\\nonumber
\widetilde{F}_3=&\,-2(D_m\mu_nu_{-1,0})(D_n\mu_mu_{-1,0})-2(D_m\mu_nv_{-1,0})(D_n\mu_mv_{-1,0})\\\nonumber
&+\eta\left\{(\Psi[u,u]-\Psi[v,v])(\Theta[u,u]-\Theta[v,v])+(\Psi[u,v]-\Psi[v,u])(\Theta[u,v]+\Theta[v,u])\right\}\\\nonumber
&-2\nu(\mu_nu_{0,0}D_nu_{0,0}+\mu_nv_{0,0}D_nv_{0,0})(\mu_nu_{-1,0}D_nu_{-1,0}+\mu_nv_{-1,0}D_nv_{-1,0}),\\\nonumber
\widetilde{G}_3=&\,\mu_m\left((D_mu_{-1,0})^2+(D_mv_{-1,0})^2\right)-\tfrac{1}2(u_{0,0}^2+v_{0,0}^2)\left\{u_{0,0}^2+v_{0,0}^2+\eta D_m^2(u_{-1,0}^2+v_{-1,0}^2)\right\},
\end{align}
and
\begin{align}\label{Theta}
\Theta[a,b]=&\,\tfrac{1}2\left\{(\mu_m\mu_na_{-1,0})D_mD_nb_{-1,0}-(D_m\mu_na_{-1,0})D_n\mu_mb_{-1,0}\right\},\\\nonumber
\Psi[a,b]=&\,2\mu_m(a_{-1,0}b_{-1,1})-{\Delta x^2}(D_ma_{-1,0})D_mb_{-1,1}.
\end{align}
None of these schemes preserves the conservation law for momentum.

In the numerical tests section, for simplicity, we restrict attention to the schemes with $\eta=\nu=0$, as these have the most compact approximation of the nonlinear term. This gives a one-parameter family $$\text{EC}_6(\alpha)=\widetilde{\mca}(\alpha\Delta_{\text{max}}^2 ,0 ,0 ),\qquad\Delta_{\text{max}} =\max{(\Delta x,\Delta t)}.$$
The scheme $\text{EC}_6(0)$ was introduced by Delfour \textit{et al.} in \cite{Delfour} by applying a Crank--Nicolson type method to a simple space discretization. In \cite{MatsuoFur}, the same scheme is obtained by using a complex discrete variational derivative method.

\subsubsection*{1-parameter family of local momentum-conserving methods}
We now derive schemes that preserve the local conservation laws for charge (\ref{NLSCL1}) and momentum (\ref{NLSCL2}). To further simplify the symbolic solution of
$$\mathsf{E}(\widetilde{\mca}\widetilde{\mcq}_1)=\mathbf{0},\qquad\mathsf{E}(\widetilde{\mca}\widetilde{\mcq}_2)=\mathbf{0},$$ 
we discretize the characteristics as follows:
$$\widetilde{\mcq}_1=(2\mu_nu_{0,0},-2\mu_nv_{0,0})^T,\qquad \widetilde{\mcq}_2=(2D_m\mu_m\mu_nv_{0,0},2D_m\mu_m\mu_nu_{0,0})^T.$$
This yields a one-parameter family of methods,
$$\widetilde{\mca}(\lambda)=(\widetilde{A}[u,v],\, \widetilde{A}[-v,u])^T,$$
where
\begin{align*}
\widetilde{A}[a,b]=&\,D_n(a_{0,0}+\lambda D_m^2a_{-1,0})+D_m^2\mu_nb_{-1,0}\\
&+\left[\mu_na_{0,0}\left(\mu_na_{0,0}+\tfrac{1}2\Delta x^2D_m^2\mu_na_{-1,0}\right)+\mu_nb_{0,0}\left(\mu_nb_{0,0}+\tfrac{1}2\Delta x^2D_m^2\mu_nb_{-1,0}\right)\right]\mu_n b_{0,0},
\end{align*}
and $\lambda=\mathcal{O}(\Delta x^2,\Delta t^2)$. The discretized conservation laws are
\begin{equation}\label{mmcl} \widetilde{\mca}\widetilde{\mcq}_1=D_m\widetilde{F}_1+D_n\widetilde{G}_1=0,\qquad\widetilde{\mca}\widetilde{\mcq}_2=D_m\widetilde{F}_2+D_n\widetilde{G}_2=0,
\end{equation}
with $\widetilde{F}_1$ and $\widetilde{G}_1$ as defined in (\ref{dfec}) and
\begin{align*}
\widetilde{F}_2=&\,(D_m\mu_nu_{-1,0})^2+(D_m\mu_nv_{-1,0})^2+(\mu_m\mu_nv_{-1,0})D_n\mu_mu_{-1,0}-(\mu_m\mu_nu_{-1,0})D_n\mu_mv_{-1,0}\\
&+\tfrac{1}2\{(\mu_nu_{0,0})\mu_nu_{-1,0}+(\mu_nv_{0,0})\mu_nv_{-1,0}\}^2\\
&+\left(\lambda-\tfrac{1}4\Delta x^2\right)\left\{(D_mD_nu_{-1,0})D_m\mu_nv_{-1,0}-(D_mD_nv_{-1,0})D_m\mu_nu_{-1,0}\right\},\\
\widetilde{G}_2=&\,u_{0,0}(D_m\mu_mv_{-1,0})\!-\!v_{0,0}(D_m\mu_mu_{-1,0})\!+\!\lambda\left\{\!(D_m^2u_{-1,0})D_m\mu_mv_{-1,0}\!-\!(D_m^2v_{-1,0})D_m\mu_mu_{-1,0}\right\}.
\end{align*}
We denote these schemes by
$$\text{MC}_6(\beta)=\widetilde{\mca}(\beta\Delta_{\text{max}}^2),\qquad\Delta_{\text{max}}=\max{(\Delta x,\Delta t)}.$$

\subsubsection*{Conservative discretizations of the Ablowitz--Ladik equation}
One can find other compact schemes that preserve two conservation laws of the NLS equation by discretizing the well-known Ablowitz--Ladik model,
\begin{equation*}
(\text{AL}[U,V],\,\text{AL}[-V,U])^T=\mathbf{0},\quad \text{AL}[a,b]=\mathrm{d}a_0/{\mathrm{d}t}+D_m^2b_{-1}+(a_0^2+b_0^2)(b_1+b_{-1})/2,
\end{equation*}
where $U_i(t)\approx u(x_i,t)$ and $V_i(t)\approx v(x_i,t)$. Therefore, we examine discretizations of the form $$\widetilde{\mca}=(\widetilde{A}[u,v],\,\widetilde{A}[-v,u]),$$ where, for second-order accuracy,
\begin{align*}
\widetilde{A}[a,b]=&\,D_na_{0,0}+D_m^2\mu_nb_{-1,0}+\!\left\{\lambda(a_{0,0}^2+a_{0,1}^2)\!+\!(1\!-\!2\lambda)a_{0,0}a_{0,1}\!+\!\theta(b_{0,0}^2+b_{0,1}^2)\!+\!(1\!-\!2\theta)b_{0,0}b_{0,1}\right\}\\\nonumber
&\times\{\phi(b_{1,1}+b_{-1,0})+(1/2-\phi)(b_{-1,1}+b_{1,0})\}.
\end{align*}
We consider the most general approximations of the characteristics, whose components are of the form (\ref{linapp}).
By solving 
\begin{equation}\label{conscondAL}
\mathsf{E}(\widetilde{\mca}\widetilde{\mcq}_1)=\mathbf{0},\qquad \mathsf{E}(\widetilde{\mca}\widetilde{\mcq}_2)=\mathbf{0},
\end{equation}
we find three different schemes. The first is
$$\text{MC-AL}=\widetilde{\mca}=(\widetilde{A}[u,v],\,\widetilde{A}[-v,u])^T=\mathbf{0}$$ with
\begin{align*}
\widetilde{A}[a,b]=&\,D_na_{0,0}+D_m^2\mu_nb_{-1,0}+\tfrac{1}2\left\{(\mu_na_{0,0})^2+(\mu_nb_{0,0})^2\right\}\mu_n(b_{-1,0}+b_{1,0}).
\end{align*}
The local conservation laws of this scheme are of the form (\ref{mmcl}) with
\begin{align*}\nonumber
\widetilde{\mcq}_1=&\left(\mu_n(u_{-1,0}+u_{1,0}),\,\mu_n(v_{-1,0}+v_{1,0})\right)^T,\\
\widetilde{F}_1=&\,2(\mu_m\mu_nu_{-1,0})D_m\mu_nv_{-1,0}-2(\mu_m\mu_nv_{-1,0})D_m\mu_nu_{-1,0}-\Delta x^2(\Theta[u,u]+\Theta[v,v]),\\
\widetilde{G}_1=&\,\tfrac{1}2\left\{u_{0,0}(u_{-1,0}+u_{1,0})+v_{0,0}(v_{-1,0}+v_{1,0})\right\},\\
\widetilde{\mcq}_2=&\left(2D_m\mu_m\mu_nu_{-1,0},\,2D_m\mu_m\mu_nv_{-1,0}\right)^T,\\
\widetilde{F}_2=&\,(D_n\mu_mu_{-1,0})\mu_n\mu_mv_{-1,0}-(D_n\mu_mv_{-1,0})\mu_m\mu_nu_{-1,0}+(D_m\mu_nu_{-1,0})^2+(D_m\mu_nv_{-1,0})^2\\
&+\tfrac{1}2\{(\mu_nu_{-1,0})^2+(\mu_nv_{-1,0})^2\}\{(\mu_nu_{0,0})^2+(\mu_nv_{0,0})^2\}\\
&+\tfrac{1}4\Delta x^2\left\{(D_m\mu_nu_{-1,0})D_mD_nv_{-1,0}-(D_m\mu_nv_{-1,0})D_mD_nu_{-1,0}\right\},\\
\widetilde{G}_2=&\,u_{0,0}D_m\mu_mv_{-1,0}-v_{0,0}D_m\mu_mu_{-1,0},
\end{align*}
where $\Theta[a,b]$ is defined in (\ref{Theta}). This scheme does not preserve the energy conservation law. 

The other two schemes that satisfy (\ref{conscondAL}) are
$$\text{M/EC-AL}(s)\equiv\widetilde{\mca}{(s)}=(\widetilde{A}[u,v],\,\widetilde{A}[-v,u])=\mathbf{0},\qquad s\in\{0,1\},$$ with
\begin{align*}
\widetilde{A}[a,b]&=D_na_{0,0}+D_m^2\mu_nb_{-1,0}+\tfrac{1}2\mu_n(a_{0,0}^2+b_{0,0}^2)(b_{-2s+1,1}+b_{2s-1,0}).
\end{align*}
For each of these schemes, there is a $\widetilde{\mcq}_3$ that satisfies
\[ \mathsf{E}(\widetilde{\mca}{(s)}\widetilde{\mcq}_3)=\mathbf{0}.\]
Therefore they have three discrete conservation laws, as follows:
\begin{align*}\nonumber
\widetilde{\mcq}_1=&\left(u_{-2s+1,1}+u_{2s-1,0},v_{-2s+1,1}+v_{2s-1,0}\right)^T,\\
\widetilde{F}_1=&\,(u_{-s,1}+u_{s-1,0})D_m\mu_nv_{-1,0}-(v_{-s,1}+v_{s-1,0})D_m\mu_nu_{-1,0}\\
&-\Delta x^2(\Theta[u,u]+\Theta[v,v])+\frac{(-1)^s}2\Delta x\Delta t\left\{(D_nu_{-1,0})D_nu_{0,0}+(D_nv_{-1,0})D_nv_{0,0}\right\},\\
\widetilde{G}_1=&\,\tfrac{1}2\left\{u_{0,0}(u_{1,0}+u_{-1,0})+v_{0,0}(v_{1,0}+v_{-1,0})\right\}\\
&+\frac{(-1)^s}2\Delta x\Delta t\left\{(D_m\mu_mu_{-1,0})D_m^2v_{-1,0}-(D_m\mu_mv_{-1,0})D_m^2u_{-1,0}\right\},\\
\widetilde{\mcq}_2=&\left(\frac{1}{\Delta x}(u_{1,1-s}-u_{-1,s}),\frac{1}{\Delta x}(v_{1,1-s}-v_{-1,s})\right)^T,\\
\widetilde{F}_2=&\,\tfrac{1}{2}\left\{(v_{-s,1}+v_{s-1,0})D_n\mu_mu_{-1,0}-(u_{-s,1}+u_{s-1,0})D_n\mu_mv_{-1,0}\right\}\\
&+(D_m\mu_nu_{-1,0})^2+(D_m\mu_nv_{-1,0})^2+\tfrac{1}2\{(\mu_nu_{-1,0})^2+(\mu_nv_{-1,0})^2\}\{(\mu_nu_{0,0})^2+(\mu_nv_{0,0})^2\}\\
&+\!\Omega_s[u,v]\!-\!\Omega_s[v,u]\!+\!\Lambda_s[u,u]\!+\!\Lambda_s[u,v]\!+\!\Lambda_s[v,u]\!+\!\Lambda_s[v,v]\!-\!\frac{\Delta t}{\Delta x}(-1)^s\!\left(\Theta[u,u]-\Theta[v,v]\right),\\
\widetilde{G}_2=&\,u_{0,0}D_m\mu_mv_{-1,0}-v_{0,0}D_m\mu_mu_{-1,0}+\frac{\Delta t}{8\Delta x}(-1)^s\left\{(u_{0,0}^2+v_{0,0}^2)(u_{-1,0}^2+u_{1,0}^2+v_{-1,0}^2+v_{1,0}^2)\right.\\
&\left.+2(u_{-1,0}+u_{1,0})D_m^2u_{-1,0}+2(v_{-1,0}+v_{1,0})D_m^2v_{-1,0}\right\},\\
\widetilde{\mcq}_3=&\left(-\frac{2}{\Delta t}(v_{-2s+1,1}-v_{2s-1,0}),-\frac{2}{\Delta t}(u_{-2s+1,1}-u_{2s-1,0})\right)^T,\\
\widetilde{F}_3=&\,-\frac{2}{\Delta t}\left\{(u_{-s,1}-u_{s-1,0})D_m\mu_nu_{-1,0}+(v_{-s,1}-v_{s-1,0})D_m\mu_nv_{-1,0}\right\}-\Phi_s[u,v]+\Phi_s[v,u]\\
&-\frac{\Delta x}{2\Delta t}(-1)^s\left\{(u_{s-1,0}^2+v_{s-1,0}^2)\mu_n(u_{-s,0}^2+v_{-s,0}^2)+(u_{-s,1}^2+v_{-s,1}^2)\mu_n(u_{s-1,0}^2+v_{s-1,0}^2)\right\},\\
\widetilde{G}_3=&\,-\tfrac{1}4(u_{0,0}^2+v_{0,0}^2)(u_{-1,0}^2+u_{1,0}^2+v_{-1,0}^2+v_{1,0}^2)+(D_mu_{-1,0})D_mu_{0,0}+(D_mv_{-1,0})D_mv_{0,0}\\
&+\frac{2\Delta x}{\Delta t}(-1)^s\left\{v_{0,0}D_m\mu_mu_{-1,0}-u_{0,0}D_m\mu_mv_{-1,0}\right\},
\end{align*}
where $\Theta[a,b]$ is defined in (\ref{Theta}) and
\begin{align*}
\Omega_s[a,b]=&\,\tfrac{1}4(-1)^s\Delta x\Delta t\left\{(D_n\mu_ma_{-1,0})D_mD_nb_{-1,0}+\tfrac{1}4\Delta x^2(D_m\mu_na_{-1,0})D_mD_nb_{-1,0}\right\},\\
\Lambda_s[a,b]=&\,\tfrac{1}8(-1)^s\Delta x\Delta t\left\{\mu_m\mu_n(a_{-1,0}^2)D_mD_n(b_{-1,0}^2)-D_n\mu_m(a_{-1,0}^2)D_m\mu_n(b_{-1,0}^2)\right\}\\
&+\tfrac{1}8\Delta t^2\left\{\mu_n(a_{-1,0}^2)(D_nb_{0,0})^2+\mu_n(b_{0,0}^2)(D_na_{-1,0})^2\right\},\\
\Phi_s[a,b]=&\,\frac{\Delta x}{\Delta t}(-1)^s\left\{b_{s-1,0}D_na_{-s,0}-a_{-s,1}D_nb_{s-1,0}\right\}.
\end{align*}
Although these schemes are second-order accurate and preserve a second-order approximation of the local charge conservation law, the local truncation error in $\widetilde{Q}_2$, $\widetilde{F}_2$ and $\widetilde{G}_2$ has terms that are $\mathcal{O}\left(\Delta t/\Delta x\right)$. On the other hand, terms that are $\mathcal{O}\left(\Delta x/\Delta t\right)$ appear in the local truncation error of $\widetilde{Q}_3$, $\widetilde{F}_3$ and $\widetilde{G}_3$.
So although these schemes have three discrete conservation laws, at most one of the conservation laws of momentum and energy converges to the continuous analogue, due to the incompatibility of the requirements $\Delta x\ll \Delta t$ and $\Delta t\ll\Delta x$, respectively. The M/EC-AL schemes with $\Delta t\ll \Delta x$ are the only local mass and energy conserving schemes that can be obtained as a one-step time discretization of the Ablowitz--Ladik system.

\subsection{Numerical Tests}

We now consider two different benchmark problems for (\ref{NLS}) subject to periodic boundary conditions, to compare the schemes developed in the previous section with other schemes from the literature. The multisymplectic scheme introduced in \cite{Chen2000}, which is equivalent to the Preissman box scheme, amounts to
\begin{equation*}
\text{MS}=(\text{MS}[u,v],\text{MS}[-v,u])=\mathbf{0},
\end{equation*}
with \begin{equation*}
\text{MS}[a,b]=D_n\mu_m^2a_{-1,0}+D_m^2\mu_nb_{-1,0}+\mu_m\left\{\left[(\mu_m\mu_na_{-1,0})^2+(\mu_m\mu_nb_{-1,0})^2\right]\mu_m\mu_nb_{-1,0}\right\}.
\end{equation*}
A backward error analysis of this method can be found in \cite{ISbea}.

In \cite{Chen2000}, the authors use a method of lines (MoL) approach to derive symplectic integrators of the NLS equation. One of these is based on the following semi-discretization of the Hamiltonian functional (\ref{Hamiltonian})
\begin{equation}\label{HamHBVM}
H(t)=\tfrac{1}{2}\Delta x\sum_i\left\{U_{i}D_m^2U_{i-1}+V_{i}D_m^2V_{i-1}-\tfrac{1}2(U_{i}^2+V_{i}^2)^2\right\},
\end{equation}
where $U_i\approx u(x_i,t)$ and $V_i\approx v(x_i,t),$ yielding the system of ODEs
\begin{equation}\label{ODEsys}
\left(\begin{array}{c}\mathrm{d}U_i/\mathrm{d}t\\\mathrm{d}V_i/\mathrm{d}t\end{array}\right)=\left(\begin{array}{cc}
0&1\\
-1&0
\end{array}\right) \left(\begin{array}{c} \mathrm{d}H/\mathrm{d}U_i \\\mathrm{d}H/\mathrm{d}V_i \end{array}\right).
\end{equation}
The second-order symplectic integrator proposed in \cite{Chen2000} is then obtained by applying the implicit midpoint rule to the semidiscretization (\ref{ODEsys}). This scheme, which we call MoL-M, is equivalent to one introduced earlier by Sanz-Serna and Verwer in \cite{S-SV86}; it conserves all quadratic invariants of (\ref{ODEsys}), including the global charge, $\Delta x\sum_i (U_{i}^2+V_{i}^2)$. Finally, we apply a Hamiltonian Boundary Value Method with $2$ stages and degree $1$, denoted HBVM($2,1$), to (\ref{ODEsys}). More generally, for any $k\geq s$, HBVM($k,s$) is a Runge-Kutta method of order $2s$ that exactly preserves polynomial Hamiltonians of degree at most $2k/s$; see \cite{BIbook,BarlettiNLS} for details on this class of methods and their application to the NLS equation. We denote the resulting scheme by MoL-HBVM(2,1); it preserves the Hamiltonian (\ref{HamHBVM}) and the global energy $\widetilde{G}_3=2H$.

For each benchmark problem, the computational cost is roughly the same for all of the schemes, so the comparisons are based mainly on the solution error at the final time $t=T$, evaluated as
\begin{equation}\label{solerrNLS}
\left.\sqrt{\frac{\|u-u_{\mathrm{exact}}\|^2+\|v-v_{\mathrm{exact}}\|^2}{\|u_{\mathrm{exact}}\|^2+\|v_{\mathrm{exact}}\|^2}}\right\vert_{t=T}.
\end{equation}
For a grid with $M$ nodes in space and $N$ nodes in time, the error in the conservation laws is
\begin{align*}
\text{Err}_\alpha=&\,\Delta x\max_{j=1,\ldots,N}\left|\sum_{i=1}^M \left(\widetilde{G}_\alpha(x_i,t_j)-\widetilde{G}_\alpha(x_i,t_1)\right)\right|,\qquad \alpha=1,2,3.
\end{align*}
Whenever any of the discrete densities $\widetilde{G}_1$, $\widetilde{G}_2$ and $\widetilde{G}_3$ are undefined by a scheme, we evaluate the corresponding error as
\begin{align*}
\text{Err}_1=&\,\Delta x\max_{j=1,\ldots,N}\left|\sum_{i=1}^M \left\{(u_{i,j}^2+v_{i,j}^2)-(u_{i,0}^2+v_{i,0}^2)\right\}\right|,\\
\text{Err}_2=&\,\Delta x\max_{j=1,\ldots,N}\left|\sum_{i=1}^M\!\left\{  u_{i,j}D_m\mu_mv_{i-1,j}\!-\!v_{i,j}D_m\mu_mu_{i-1,j}\!-\!u_{i,0}D_m\mu_mv_{i-1,0}\!+\!v_{i,0}D_m\mu_mu_{i-1,0}\right\}\right|,\\
\text{Err}_3=&\,\Delta x\max_{j=1,\ldots,N}\left|\sum_{i=1}^M \left\{ \mu_m\left((D_mu_{i-1,j})^2+(D_mv_{i-1,j})^2\right)-\tfrac{1}2(u_{i,j}^2+v_{i,j}^2)^2\right.\right.\\
&\left.\qquad\qquad\quad-\mu_m\left((D_mu_{i-1,0})^2+(D_mv_{i-1,0})^2\right)+\tfrac{1}2(u_{i,0}^2+v_{i,0}^2)^2\right\}\Bigg|,
\end{align*}
with $u_{i,j}\simeq u(a+i\Delta x, j\Delta t)$.

The first benchmark problem is defined on $\Omega=[-20,20]\times [0,T]$ with the initial condition
\[\psi(x,0)=\sqrt{2}\,\sech{(x-d)}\exp\{\mathrm{i}c(x-d)\}. \]
The exact solution is
\[\psi(x,t)=\sqrt{2}\,\sech(x-d-2ct)\exp\{\mathrm{i}c(x-d)-\mathrm{i}(c^2-1)t)\}, \]
and $|\psi(x,t)|$ is a single soliton.
We set $d=-5$ and use step lengths $\Delta x=0.1$ and $\Delta t=0.02$.

To see how the accuracy varies over a long time, we choose $c=0.05$ and $T=100$ and evaluate the solution error at $t=5,10,\ldots,100.$ The plot on the top-left of Figure~\ref{errorgrownls} shows that for EC$_6(0)$ and MC$_6(0)$, this error accumulates linearly with respect to $t$, as it is to be expected when the numerical method preserves invariants \cite{defrutosSS,DuranSS}. The error given by each of M/EC-AL(0), M/EC-AL(1) and MC-AL also grows linearly, at a similar rate to that of MC$_6(0)$. By contrast, the error of the energy conserving method EC$_6(0)$ accumulates far more slowly. The remaining plots in  Figure~\ref{errorgrownls} show the error in the invariants of EC$_6(0)$, MC$_6(0)$, and M/EC-AL(1) at $t=0,0.1,\ldots,100$. Among the five schemes considered, only EC$_6(0)$ and MC-AL do not show a visible drift in the error of any of the conserved invariants.

\begin{figure}[ht!]
\begin{center}
	\definecolor{mycolor1}{rgb}{0.49020,0.18039,0.56078}%
\definecolor{mycolor2}{rgb}{0.00000,0.49804,0.00000}%

\begin{tikzpicture}

\begin{axis}[%
width=2in,
height=1.5in,
at={(0.758in,0.481in)},
scale only axis,
xmin=0,
xmax=100,
xlabel style={font=\color{white!15!black}},
xlabel={$t$},
tick label style={/pgf/number format/fixed},
scaled y ticks=false,
ymin=0,
ymax=0.8,
title={Solution error},
axis background/.style={fill=white},
legend style={legend cell align=left, align=left, draw=white!15!black,at={(0.50,0.95)}}
]
\addplot [color=blue, line width=1.0pt, mark size=3.0pt, mark=o, mark options={solid, blue}]
  table[row sep=crcr]{%
5	0.0100789182439332\\
10	0.0197912387828924\\
15	0.0290795575823996\\
20	0.0384486842740348\\
25	0.0482733935391748\\
30	0.0582864910839237\\
35	0.0680406797824304\\
40	0.0775731247632275\\
45	0.0869561901397687\\
50	0.0965490329381769\\
55	0.106499400790895\\
60	0.116259934647979\\
65	0.125716445609884\\
70	0.13526043524037\\
75	0.145216423989688\\
80	0.155084489718308\\
85	0.164529548306453\\
90	0.173911822036089\\
95	0.183324921195944\\
100	0.19298961302975\\
};
\addlegendentry{EC$_6(0)$}

\addplot [color=red, line width=1.0pt, mark size=2.0pt, mark=diamond, mark options={solid, red}]
  table[row sep=crcr]{%
5	0.0308363164369454\\
10	0.0609649721639928\\
15	0.0902899857650758\\
20	0.119866784105851\\
25	0.150215374822164\\
30	0.180879985754904\\
35	0.210930619105112\\
40	0.240444362497806\\
45	0.269954403890746\\
50	0.299753726946062\\
55	0.329909201324824\\
60	0.359526157980276\\
65	0.388929933120999\\
70	0.418358809500871\\
75	0.448371980607553\\
80	0.477934034596842\\
85	0.506669250532288\\
90	0.535108865067328\\
95	0.563851832278183\\
100	0.593052330330811\\
};
\addlegendentry{MC$_6(0)$}

\end{axis}

\end{tikzpicture}%
	\input{EC6NLSclawerr.tex}
	\input{MC6NLSclawerr.tex}
	\input{MEC1clawerr.tex}
\end{center}
\caption{Solution error and invariants error of different schemes.}\label{errorgrownls}
\end{figure}

To illustrate the parameter dependence, we now set $c=2.5$ and $T=2$; the optimal values of the free parameters for $\text{EC}_6(\alpha)$ and $\text{MC}_6(\beta)$ are $\alpha=0.132$ and $\beta=0.043$, respectively.

Table~\ref{comp1solnls} shows that all the schemes introduced in the previous section preserve two conservation laws. The M/EC-AL schemes preserve three invariants, but as $\Delta x>\Delta t$, we do not expect the discrete energy to be a good approximation of the continuous one. 

\begin{figure}[htbp]
\begin{center}
	\input{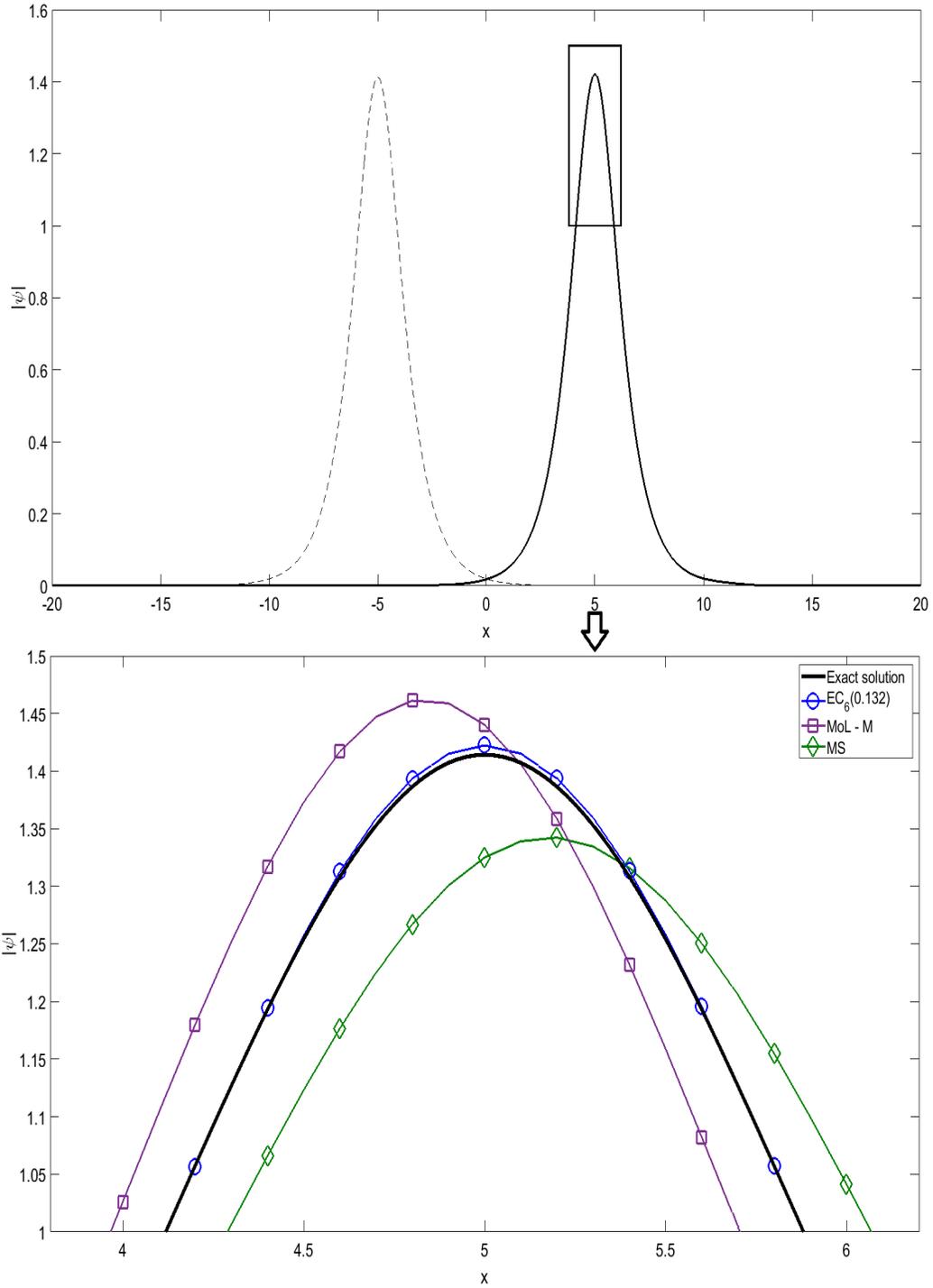}
\end{center}
\caption{Single soliton problem for the NLS equation. Top: Modulus of the initial condition (dashed line) and of the solution of $\mbox{EC}_6(0.132)$ with $\Delta x=0.1$ and $\Delta t=0.02$ over $[-20,20]$ at time $T=2$ (solid line). Bottom: Top of the soliton at time $T=2$. Modulus of the reference solution (bold line) and the solutions of $\mbox{EC}_6(0.132)$ (circles), MoL-M (squares), and MS (diamonds). Markers are at every second point.}\label{ssNLS}
\end{figure}

All of the new schemes compare well with existing methods; some do much better, the most accurate being EC$_6$(0.132).
The upper part of Figure~\ref{ssNLS} shows the modulus of the initial condition and the numerical solution given by EC$_6$(0.132) at time $T=2$. The lower plot shows the exact solution and the numerical solutions given by EC$_6$(0.132), MoL-M and MS, close to the top of the soliton. We do not show the solutions from MoL-HBVM(2,1) and EC$_6$(0), as they almost overlap the one given by MoL-M. The solution of EC$_6$(0.132) is by far the closest to the exact solution and matches both the amplitude and phase of the soliton well.

\begin{table}[ht]
\caption{Single soliton problem for the NLS equation; $\Delta x=0.1$ and $\Delta t=0.02$.}
\label{comp1solnls}
\small
\centerline{\begin{tabular}{|c|c|c|c|c|c|}
\hline
Method &  $\text{Err}_1$ & $\text{Err}_2$ & $\text{Err}_3$  &  Error in the solution\\ 
\hline
Delfour \textit{et al.}; $\mbox{EC}_6(0)$	& 1.11e-14   & 0.0032 & 7.11e-14 &  0.1887\\
\hline
$\mbox{EC}_{6}(0.132)$	& 1.07e-14  & 1.86e-04 & 7.46e-14 &  0.0083\\	
\hline
$\mbox{MC}_6(0)$	& 9.33e-15 & 1.42e-14 & 0.0016 & 0.0952\\
\hline
$\mbox{MC}_6(0.043)$	& 8.88e-15  & 2.31e-14 & 4.62e-04 & 0.0741 \\
\hline
$\mbox{M/EC-AL}(1)$	& 8.44e-15  & 2.13e-14 & 2.13e-13 & 0.0729 \\	
\hline
$\mbox{M/EC-AL}(0)$	& 9.77e-15  & 1.42e-14 & 1.42e-13 & 0.1267 \\	
\hline
$\mbox{MC-AL}$	& 1.07e-14  & 2.13e-14 & 8.92e-05 & 0.0765 \\	
\hline
MS & 6.63e-04  & 0.0016 & 0.0077 &  0.3023\\
\hline
MoL - M	& 1.11e-14  & 0.0025 & 0.0023 & 0.1785 \\	
\hline
MoL - HBVM(2,1)	& 9.93e-05 & 0.0032 & 7.11e-14 & 0.1817 \\	
\hline
\end{tabular}}
\end{table}

\begin{figure}[tbp]
{\includegraphics[width=16.5cm,height=9.75cm]{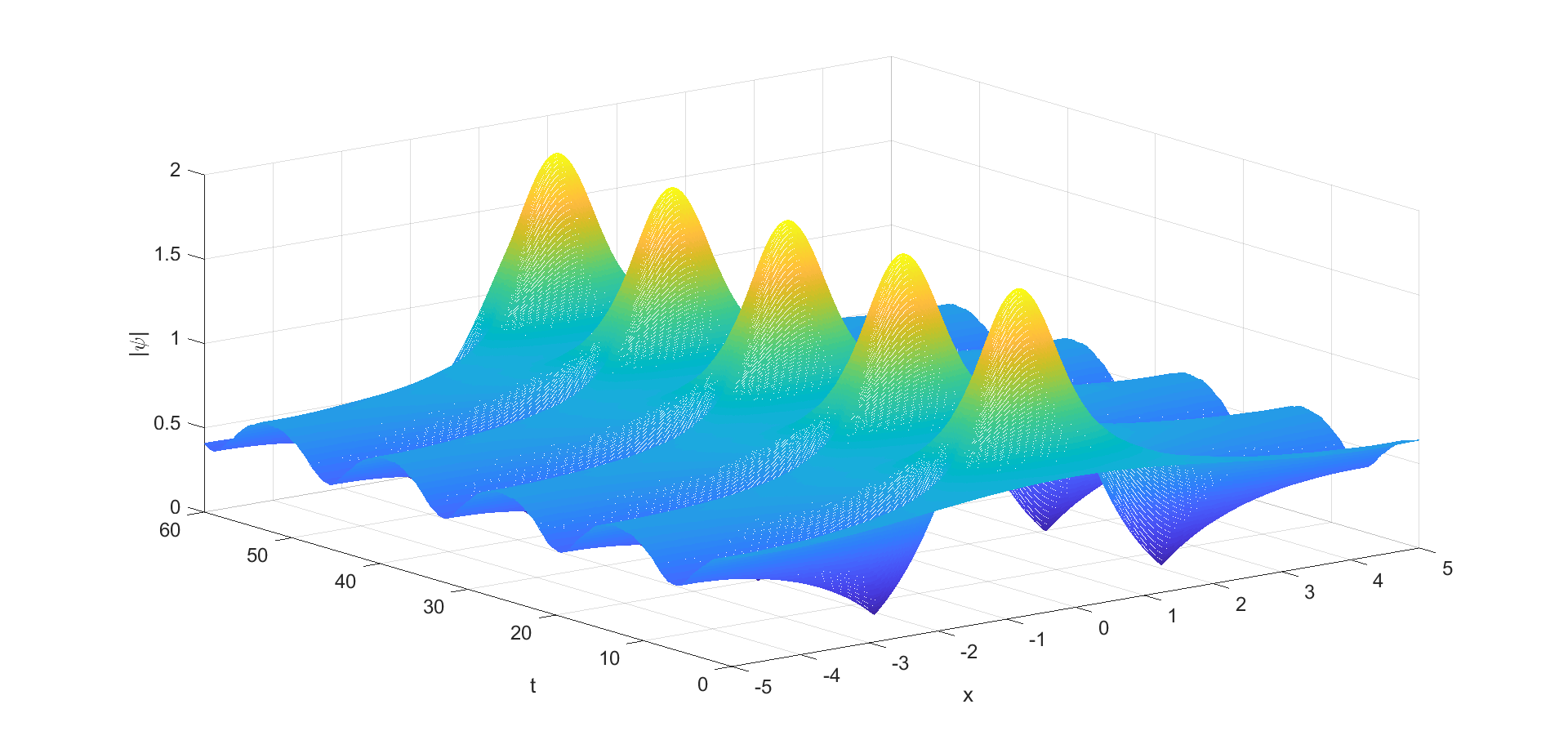}}\vspace{-0.4cm}
\caption{NLS with initial condition (\ref{Brinit}): modulus of the reference solution.  }\label{bNLS}
\end{figure}

The second benchmark problem is a breather solution on $(x,t)\in[-5,5]\times[0,60]$, with initial condition \cite{IslasSchober}
\begin{equation}\label{Brinit}
\psi(x,0)=\frac{1}{\sqrt{2}}\left\{1+0.1\cos\left(\frac{x}{\sqrt{2}}\right)\right\}.
\end{equation}
The reference solution is obtained by applying the sixth-order energy-conserving method HBVM(6,3), with $\Delta t=0.05$, to the system of ODEs (\ref{ODEsys}) resulting from a spatial discretization with step $\Delta x=0.0125$. The modulus of the solution obtained is shown in Figure~\ref{bNLS} and exhibits the expected quasi-periodic (in time) breather motion (see \cite{IslasSchober}). We compare this with the various second-order schemes defined on the much coarser grid with $\Delta x=0.05$ and $\Delta t=0.2$. On this grid, the values $\alpha=0.052$ and $\beta=0.371$ minimize the solution error of EC$_6(\alpha)$ and MC$_6(\beta)$.

\begin{table}[t]
\caption{NLS equation with initial condition (\ref{Brinit}); $\Delta x=0.05$ and $\Delta t=0.2$.}
\label{tabbreathnls}
\small
\centerline{\begin{tabular}{|c|c|c|c|c|c|}
\hline
Method &  $\text{Err}_1$ & $\text{Err}_2$ & $\text{Err}_3$  &  Error in the solution\\ 
\hline
Delfour \textit{et al.}; $\mbox{EC}_6(0)$	&  8.88e-15  & 3.43e-15 & 1.33e-14 &  0.0765\\
\hline
$\mbox{EC}_{6}(0.052)$	& 1.15e-14  & 4.83e-15 & 1.29e-14 &  0.0231\\
\hline   
$\mbox{MC}_6(0)$	& 1.07e-14 & 9.90e-15 & 0.0471 & 0.6121\\ 
\hline
$\mbox{MC}_6(0.371)$	&  7.99e-15 & 7.04e-15 & 0.0515 & 0.0966 \\
\hline
$\mbox{M/EC-AL}(1)$	& 1.24e-14  & 1.42e-14 & 1.18e-14 & 0.1062 \\	
\hline
$\mbox{M/EC-AL}(0)$	& 7.99e-15  & 2.13e-14 & 1.11e-14 & 0.1062 \\	
\hline
$\mbox{MC-AL}$	& 8.88e-15  & 1.01e-14 & 0.0523 & 0.5902 \\
\hline
MS & 0.0023  & 2.77e-15 & 0.0498 & 0.5924 \\	
\hline
MoL - M	& 8.88e-15  & 9.71e-15 & 0.0424 & 0.5806 \\	
\hline
MoL - HBVM(2,1)	& 0.0259 & 1.65e-14 & 6.84e-14 & 0.5123 \\
\hline
\end{tabular}}
\end{table}

Table~\ref{tabbreathnls} shows that every scheme preserves the momentum to machine accuracy. This is a consequence of the spatial symmetry of the schemes and the solution. Again, M/EC-AL schemes preserve all three discrete invariants (up to rounding errors), but as $\Delta t>\Delta x$, we cannot expect the discrete momentum to approximate the continuous momentum accurately.

The accuracy of all the schemes discussed in the previous section and listed in Table~\ref{tabbreathnls} is at least comparable to existing methods. The local energy-conserving schemes are particularly accurate, the best being EC$_6$(0.052).

\begin{figure}[htbp]
\begin{center}
\input{breather.tex}
\input{phase.tex}
\end{center}
\caption{NLS solution with initial condition (\ref{Brinit}) at time $T=60$, with $\Delta x=0.05$ and $\Delta t=0.2$. Modulus (top) and phase (bottom) of reference solution (bold line) and solutions of $\mbox{EC}_6(0.052)$ (circles), MoL-M (squares), and MoL-HBVM(2,1) (diamonds). Markers are at every eighth point. }\label{bcNLS}
\end{figure}

Figure~\ref{bcNLS} shows the modulus (top) and phase, $\theta=\arctan{(v/u)}$, (bottom) of the reference solution and the numerical solutions given by EC$_6$(0.052), MoL-M and MoL-HBVM(2,1). The solution of MS is not shown, as it almost overlaps the solution of MoL-M. The new scheme reproduces both the oscillations of the breather motion and the phase of the complex-valued solution particularly well.

\begin{figure}[htbp]
\begin{center}
\definecolor{mycolor1}{rgb}{0.00000,0.50196,0.00000}%
\definecolor{mycolor2}{rgb}{0.49020,0.18039,0.56078}%
\begin{tikzpicture}
\begin{axis}[%
width=5in,
height=3in,
at={(1in,3.2in)},
scale only axis,
xmin=-0.6,
xmax=0.6,
xlabel style={font=\color{white!15!black}},
xlabel={$x$},
ymin=1.38,
ymax=1.58,
ylabel style={font=\color{white!15!black}},
ylabel={$|\psi|$},
axis background/.style={fill=white},
legend style={legend cell align=left, align=left, draw=white!15!black}
]
\addplot [color=blue, line width=2.0pt, mark size=5.0pt, mark=o, mark options={solid, blue, line width=1.0pt}, mark repeat = 3]
  table[row sep=crcr]{%
-0.9	0.988831822131741\\
-0.85	1.0305390025756\\
-0.8	1.0724518700605\\
-0.75	1.1143497332634\\
-0.7	1.15598235288547\\
-0.65	1.19705229651427\\
-0.6	1.23720926314899\\
-0.55	1.27605992321444\\
-0.5	1.31319021004172\\
-0.45	1.34819200213781\\
-0.4	1.38068485303419\\
-0.35	1.41032622711947\\
-0.3	1.43680917165235\\
-0.25	1.45985200840412\\
-0.2	1.47918801382226\\
-0.15	1.49456279477871\\
-0.0999999999999996	1.50574337788879\\
-0.0499999999999998	1.51253757595809\\
0	1.51481733647732\\
0.0499999999999998	1.51253757595819\\
0.0999999999999996	1.505743377889\\
0.15	1.49456279477901\\
0.2	1.47918801382266\\
0.25	1.45985200840461\\
0.3	1.43680917165292\\
0.35	1.41032622712012\\
0.4	1.38068485303491\\
0.45	1.34819200213858\\
0.5	1.31319021004255\\
0.55	1.27605992321531\\
0.6	1.23720926314989\\
0.65	1.1970522965152\\
0.7	1.1559823528864\\
0.75	1.11434973326435\\
0.8	1.07245187006145\\
0.85	1.03053900257655\\
0.9	0.988831822132682\\
};
\addlegendentry{EC$_{6}(0.052)$}

\addplot [color=mycolor1, line width=2.0pt, mark size=5.0pt, mark=x, mark options={solid, mycolor1, line width=1.0pt}, mark repeat = 3]
  table[row sep=crcr]{%
-0.9	0.986078225902971\\
-0.85	1.02926812698962\\
-0.8	1.07278923862049\\
-0.75	1.11642174690406\\
-0.7	1.15991035433181\\
-0.65	1.2029459017718\\
-0.6	1.24515932900526\\
-0.55	1.28613156669759\\
-0.5	1.32541628203067\\
-0.45	1.36256736989844\\
-0.4	1.39716179452837\\
-0.35	1.42881120054083\\
-0.3	1.45716119770774\\
-0.25	1.48188287837454\\
-0.2	1.50266449491358\\
-0.15	1.51921092828499\\
-0.0999999999999996	1.5312548460826\\
-0.0499999999999998	1.53857795468938\\
0	1.54103587162122\\
0.0499999999999998	1.53857795468937\\
0.0999999999999996	1.53125484608259\\
0.15	1.51921092828497\\
0.2	1.50266449491356\\
0.25	1.48188287837452\\
0.3	1.45716119770771\\
0.35	1.4288112005408\\
0.4	1.39716179452834\\
0.45	1.3625673698984\\
0.5	1.32541628203062\\
0.55	1.28613156669755\\
0.6	1.24515932900522\\
0.65	1.20294590177176\\
0.7	1.15991035433176\\
0.75	1.11642174690402\\
0.8	1.07278923862045\\
0.85	1.02926812698957\\
0.9	0.986078225902928\\
};
\addlegendentry{EC$_6(0)$}

\addplot [color=black, line width=2.0pt]
  table[row sep=crcr]{%
-0.9	0.984443282756461\\
-0.8875	0.995145270313255\\
-0.875	1.00588996031897\\
-0.8625	1.0166777639406\\
-0.85	1.0275035921071\\
-0.8375	1.03835706096613\\
-0.825	1.04922559029219\\
-0.8125	1.06009869921558\\
-0.8	1.07097109857963\\
-0.7875	1.0818428777781\\
-0.775	1.09271676385543\\
-0.7625	1.10359409580718\\
-0.75	1.11447183716545\\
-0.7375	1.1253422719217\\
-0.725	1.13619542264375\\
-0.7125	1.14702265547727\\
-0.7	1.15781930603215\\
-0.6875	1.1685848073575\\
-0.675	1.17932030528254\\
-0.6625	1.19002520490261\\
-0.65	1.20069466154715\\
-0.6375	1.21131941535066\\
-0.625	1.22188797166333\\
-0.6125	1.23238978677788\\
-0.6	1.24281759943015\\
-0.5875	1.25316761668826\\
-0.575	1.26343754720274\\
-0.5625	1.27362371062607\\
-0.55	1.28371893664926\\
-0.5375	1.29371245879662\\
-0.525	1.30359183941815\\
-0.5125	1.31334582308229\\
-0.5	1.32296654055878\\
-0.4875	1.33244990399648\\
-0.475	1.34179406388129\\
-0.4625	1.35099683695727\\
-0.45	1.36005349423846\\
-0.4375	1.36895600099357\\
-0.425	1.37769395071529\\
-0.4125	1.38625654026467\\
-0.4	1.39463447911\\
-0.3875	1.40282091040635\\
-0.375	1.41081107436129\\
-0.3625	1.41860115335578\\
-0.35	1.42618710440954\\
-0.3375	1.43356414465256\\
-0.325	1.44072705048685\\
-0.3125	1.44767089491083\\
-0.3	1.45439159731453\\
-0.2875	1.46088581793656\\
-0.275	1.46715017138139\\
-0.2625	1.47318018045463\\
-0.25	1.47896958024009\\
-0.2375	1.48451041478883\\
-0.225	1.48979395572164\\
-0.2125	1.49481205201938\\
-0.2	1.49955831769234\\
-0.1875	1.50402866762639\\
-0.175	1.50822104416549\\
-0.1625	1.51213455052048\\
-0.15	1.51576843431301\\
-0.1375	1.5191213562593\\
-0.125	1.52219117485782\\
-0.1125	1.524975209139\\
-0.0999999999999996	1.52747074570129\\
-0.0875000000000004	1.52967550870278\\
-0.0750000000000002	1.53158790123065\\
-0.0625	1.53320697885883\\
-0.0499999999999998	1.53453224516149\\
-0.0374999999999996	1.53556341067196\\
-0.0250000000000004	1.53630022923285\\
-0.0125000000000002	1.53674245455738\\
0	1.53688988864821\\
0.0125000000000002	1.53674245455847\\
0.0250000000000004	1.53630022923168\\
0.0374999999999996	1.53556341067288\\
0.0499999999999998	1.53453224516244\\
0.0625	1.53320697885866\\
0.0750000000000002	1.53158790123215\\
0.0875000000000004	1.52967550870248\\
0.0999999999999996	1.52747074570297\\
0.1125	1.52497520914008\\
0.125	1.52219117485801\\
0.1375	1.51912135626167\\
0.15	1.51576843431363\\
0.1625	1.51213455052219\\
0.175	1.50822104416827\\
0.1875	1.50402866762605\\
0.2	1.49955831769517\\
0.2125	1.49481205202122\\
0.225	1.48979395572291\\
0.2375	1.48451041479223\\
0.25	1.47896958024038\\
0.2625	1.47318018045796\\
0.275	1.46715017138383\\
0.2875	1.4608858179385\\
0.3	1.45439159731693\\
0.3125	1.4476708949139\\
0.325	1.44072705048919\\
0.3375	1.43356414465586\\
0.35	1.42618710441212\\
0.3625	1.41860115335817\\
0.375	1.41081107436534\\
0.3875	1.40282091040833\\
0.4	1.39463447911409\\
0.4125	1.38625654026736\\
0.425	1.37769395071834\\
0.4375	1.3689560009983\\
0.45	1.36005349424012\\
0.4625	1.35099683696208\\
0.475	1.34179406388382\\
0.4875	1.33244990400026\\
0.5	1.32296654056203\\
0.5125	1.31334582308668\\
0.525	1.30359183942099\\
0.5375	1.29371245880113\\
0.55	1.28371893665333\\
0.5625	1.27362371062896\\
0.575	1.26343754720652\\
0.5875	1.25316761669206\\
0.6	1.24281759943453\\
0.6125	1.23238978678296\\
0.625	1.22188797166513\\
0.6375	1.21131941535552\\
0.65	1.20069466155138\\
0.6625	1.19002520490657\\
0.675	1.17932030528649\\
0.6875	1.16858480736157\\
0.7	1.15781930603647\\
0.7125	1.14702265548101\\
0.725	1.13619542264835\\
0.7375	1.12534227192415\\
0.75	1.11447183717177\\
0.7625	1.10359409581021\\
0.775	1.09271676385938\\
0.7875	1.08184287778246\\
0.8	1.07097109858392\\
0.8125	1.06009869921958\\
0.825	1.04922559029521\\
0.8375	1.03835706097092\\
0.85	1.02750359211212\\
0.8625	1.01667776394373\\
0.875	1.00588996032362\\
0.8875	0.995145270316829\\
0.9	0.984443282760312\\
};
\addlegendentry{Reference solution}

\end{axis}

\end{tikzpicture}%
\input{phaseEC.tex}
\end{center}
\caption{NLS solution with initial condition (\ref{Brinit}) at time $T=60$, with $\Delta x=0.05$ and $\Delta t=0.2$. Modulus (top) and phase (bottom) of reference solution (bold line) and solutions of $\mbox{EC}_6(0.052)$ (circles) and $\mbox{EC}_6(0)$ (crosses). Markers are at every third (top) or eighth (bottom) point. }\label{bzNLS}
\end{figure}

\begin{figure}[htbp]
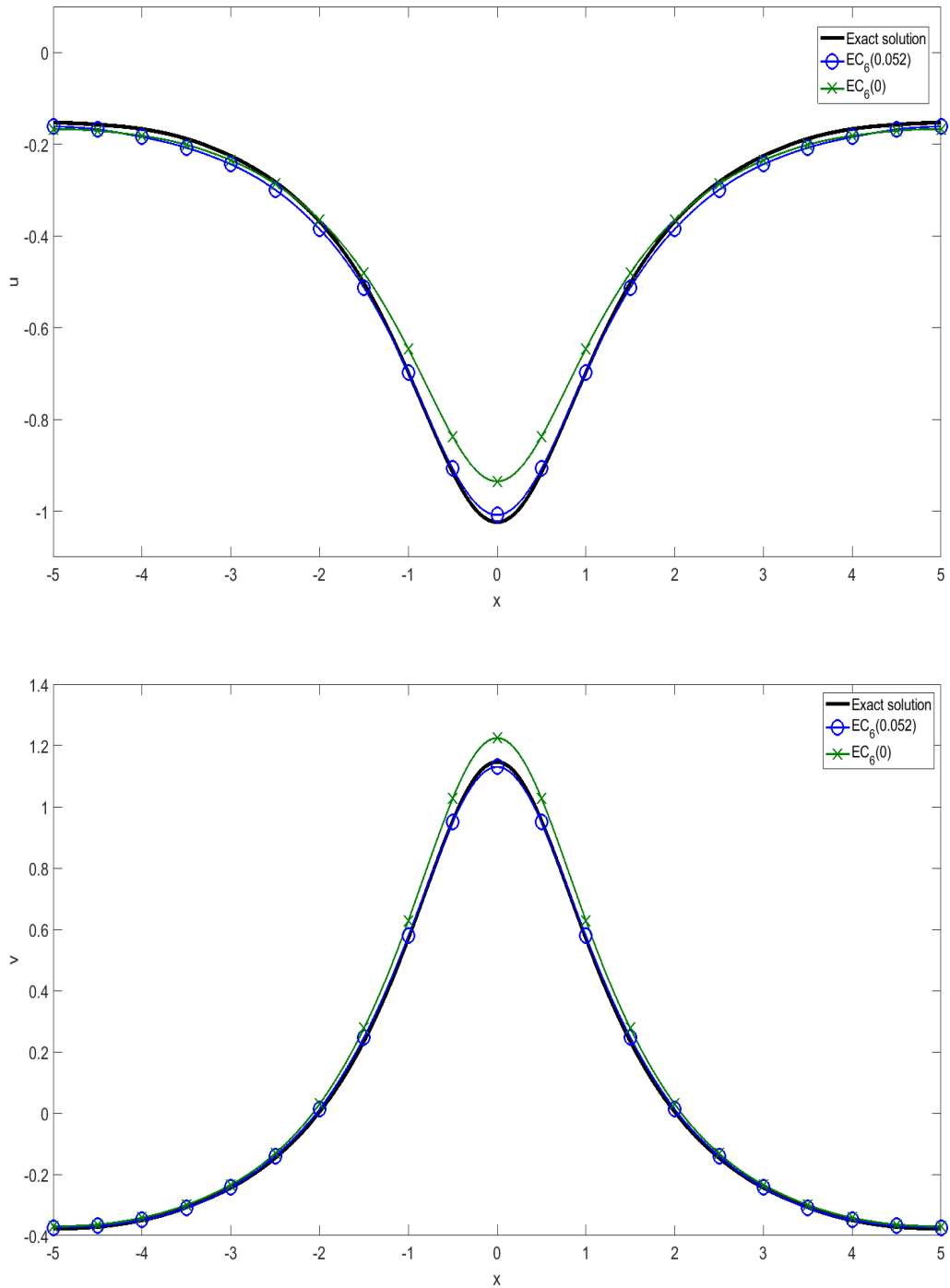

\begin{center}
\input{ubreather.tex}
\input{vbreather.tex}
\end{center}
\caption{NLS equation with initial condition (\ref{Brinit}) at time $T=60$, with $\Delta x=0.05$ and $\Delta t=0.2$. Real part (top) and imaginary part (bottom) of the reference solution (bold line), $\mbox{EC}_6(0.052)$ (circles) and $\mbox{EC}_6(0)$ (crosses). Markers are at every tenth point. }\label{buvNLS}
\end{figure}

The upper plot in Figure~\ref{bzNLS} compares the modulus of the exact solution and the numerical solutions given by EC$_6$(0.052) and by the Delfour \textit{et al.} method (EC$_6$(0)) around the maximum, where the difference between the two solutions is largest. Note that EC$_6$(0) reproduces the modulus of the solution for $x\in[-0.6,0.6]$ (the error is $\approx 0.0041$) better than EC$_6$(0.052) does (error $\approx$ 0.0221). Nevertheless, for $x\in[-0.6,0.6]$, the phase error of EC$_6(0)$ ($\approx$ 0.0818) is much larger than that of EC$_6(0.052)$ ($\approx$ 0.0025). This can be seen in the lower plot in Figure~\ref{bzNLS}. Consequently, in accordance with the results in Table~\ref{tabbreathnls}, EC$_6$(0.052) approximates the real and imaginary parts of the solution best, as shown in Figure~\ref{buvNLS}.

\section{Conclusions}\label{Conclusions}
The strategy introduced in \cite{FCHydon} for efficiently deriving bespoke finite difference methods that preserve conservation laws has been applied, for the first time, to a PDE not in Kovalevskaya form (the BBM equation) and to a system of PDEs (the NLS equation in real form). New parametrized families of conservative numerical schemes have been found for each of these equations. For the NLS equation, we have also derived new time integrators of the Ablowitz--Ladik model that preserve multiple conservation laws. 

For various benchmark problems, we have found members in each new family of schemes that give very accurate solutions compared to other schemes from the literature. More generally, free parameters in families of bespoke schemes can be chosen to reduce the error in other geometric features, including further conservation laws (see \cite{FCHydon, FCHmkdv}).

There are several ways to extend the approach described in this paper. One is to discretize in space only (reducing compactness restrictions) and combine this with known conservative time integrators. We are currently working to develop high-order schemes this way. We are also investigating generalizations to conservation laws of non-polynomial PDEs.


\begin{thebibliography}{99}
\setlength{\itemsep}{0.0cm}\setlength{\parsep}{0cm}

\bibitem{AblowitzLadik} M.\,J.\,Ablowitz, J.\,F.\,Ladik. A nonlinear difference scheme and inverse scattering. {\em Stud. Appl. Math.}, {\bf 55} (1976), 213--229.

\bibitem{An1992} I.\,M.\,Anderson, Introduction to the Variational Bicomplex, in Mathematical Aspects of Classical Field Theory, {\it Contemp. Math.}, {\bf 132} (1992), 51--73.

\bibitem{BarlettiNLS} L.\,Barletti, L.\,Brugnano, G.\,Frasca-Caccia,\,F.\,Iavernaro. Energy-conserving methods for the nonlinear Schr\"odinger equation. {\em Appl. Math. Comput.}, {\bf 318} (2018), 3--18.

\bibitem{BonaPritchardScott} J.\,L.\,Bona, W.\,G.\,Pritchard, L.\,R.\,Scott. A comparison of solutions of two model equations for long waves. {\em Lectures in Applied Mathematics}, {\bf 20} (1983), 235--267.

\bibitem{Bridges} T.\,J.\,Bridges, S.\,Reich. Multi-symplectic integrators: numerical schemes for Hamiltonian PDEs that conserve symplecticity. {\em Phys. Lett. A}, {\bf 284} (2001), 184--193.

\bibitem{BridgesReich} T.\,J.\,Bridges, S.\,Reich. Numerical methods for Hamiltonian PDEs. {\em J. Phys. A: Math. Gen.}, {\bf 39}, (2006), 5287--5320.

\bibitem{BFI15} 
L.\,Brugnano, G.\,Frasca-Caccia, F.\,Iavernaro. Energy conservation issues in the numerical solution of the semilinear wave equation.
{\em Appl. Math. Comput.}, {\bf 270} (2015), 842--870.

\bibitem{BIbook} L.\,Brugnano, F.\,Iavernaro. {\em Line Integral Methods for Conservative Problems.} Monograph and Research Notes in Mathematics. Boca Raton, FL: CRC Press, 2016.

\bibitem{Celledoni} E.\,Celledoni, V.\,Grimm, R.\,I.\,McLachlan, D.\,I.\,McLaren, D.\,O'Neal, B.\,Owren, G.\,R.\,W.\, Quispel. Preserving energy resp. dissipation in numerical PDEs using the ``Average Vector Field'' method.  {\em J. Comput. Phys.}, {\bf 231} (2012), 6770--6789.

\bibitem{Cell2} E.\,Celledoni, R.\,I.\,McLachlan, B.\,Owren, G.\,R.\,W.\,Quispel. Energy-preserving integrators and the structure of B-series. {\em Found. Comput. Math.},
{\bf 10} (2010), 673–693.

\bibitem{ChenQin} J.\,B.\,Chen, M.\,Z.\,Qin. Multi-symplectic Fourier pseudospectral method for the nonlinear Schr\"odinger equation. {\em Electron. Trans. Numer. Anal.,} {\bf 1} (2001), 193--204.

\bibitem{Chen2000} J.\,B.\,Chen, M.\,Z.\,Qin, Y.\,F.\,Tang. Symplectic and multi-symplectic methods for the nonlinear Schr\"odinger equation. {\em Comput. Math. Appl.,} {\bf 43} (2002), 1095--1106.

\bibitem{Dahlby} M.\,Dahlby, B.\,Owren. A general framework for deriving integral preserving numerical methods for
PDEs. {\em SIAM J. Sci. Comput.}, {\bf 33} (2011), 2318--2340.

\bibitem{defrutosSS} J.\,de Frutos, J.\,M.\,Sanz-Serna. Accuracy and conservation properties in numerical integration: the case of the Korteweg-de Vries equation. {\em Numer. Math.}, {\bf 75} (1997), 421--445.

\bibitem{DuranLM} A.\,Dur\'{a}n, M.\,A.\,L\'{o}pez-Marcos. Conservative numerical methods for solitary wave interactions. {\em J. Phys. A: Math. Gen.}, {\bf 36} (2003), 7761--7770.

\bibitem{DuranSS98} A.\,Dur\'{a}n, J.\,M.\,Sanz-Serna. The numerical integration of relative equilibrium solutions. Geometric theory. {\em Nonlinearity}, {\bf 11} (1998), 1547--1567.

\bibitem{DuranSS} A.\,Dur\'{a}n, J.\,M.\,Sanz-Serna. The numerical integration of relative equilibrium solutions. The nonlinear Schr\"odinger equation. {\em IMA J. Numer. Anal.}, {\bf 20} (2000), 235--261.

\bibitem{Delfour} M.\,Delfour, M.\,Fortin, G.\,Payre. Finite-difference solution of a nonlinear Schr\"odinger equation. {\em J. Comput. Phys.}, {\bf 44} (1981), 277--288.

\bibitem{DuzhinTsujishita} S.\,V.\,Duzhin, T.\,Tsujishita. Conservation laws of the BBM equation. {\em J. Phys. A: Math. Gen.}, {\bf 17} (1984), 3267--3276.

\bibitem{FT87} L.\,Faddeev, L.\,Takhtajan. {\em Hamiltonian Methods in the Theory of Solitons}, 1st ed., Springer-Verlag, Berlin, 1987.

\bibitem{FC18} G.\,Frasca-Caccia. Bespoke finite difference methods that preserve two local conservation laws of the modified KdV equation. {\em AIP Conf. Proc.}, {\bf 2116}
(2019), 140004.

\bibitem{FCHydon} G.\,Frasca-Caccia, P.\,E.\,Hydon. Simple bespoke preservation of two conservation laws. {\em IMA J. Numer. Anal.}, {\bf 40} (2020), 1294-1329.

\bibitem{FCHmkdv} G.\,Frasca-Caccia, P.\,E.\,Hydon. Locally conservative finite difference schemes
for the modified KdV equation. {\em J. Comput. Dyn.}, {\bf 6} (2019), 307-323.

\bibitem{Gonz} O.\,Gonzales. Time integration and discrete Hamiltonian systems. {\em J. Nonlinear Sci.}, {\bf 6} (1996), 449–467.

\bibitem{Grant} 
{T.\,J. Grant. } 
Bespoke finite difference schemes that preserve multiple conservation laws.
{\em LMS J. Comput. Math.}, {\bf 18} (2015), 372--403.

\bibitem{GrantHydon} 
{T.\,J.\,Grant,  P.\,E.\, Hydon. } 
 Characteristics of conservation laws for difference equations.
 {\em Found. Comput. Math.}, {\bf 13} (2013), 667--692.

\bibitem{Heitz} C.\,Heitzinger, C.\,Ringhofer, S.\,Selberherr. Finite difference solutions of the nonlinear Schr\"odinger equation and their conservation of physical quantities. {\em Commun. Math. Sci.}, {\bf 5} (2007), 779--788.

\bibitem{Hydon2001} P.\,E.\,Hydon. Conservation laws of partial difference equations with two independent variables.
{\em J. Phys. A}, {\bf 34} (2001), 10347--10355.

\bibitem{Hydonbook} P.\,E.\,Hydon. {\em Difference Equations by Differential Equation Methods}. Cambridge University Press, Cambridge, 2014.

\bibitem{HydonMans}  P.\,E.\,Hydon, E.\,L.\,Mansfield.  A variational complex for difference equations. {\em Found. Comput. Math.}, {\bf 4} (2004), 187--217.

\bibitem{IslasSchober} A.\,L.\,Islas, C.\,M.\,Schober. On the preservation of phase space structure under multisymplectic discretization. {\em J. Comput. Phys.}, {\bf 197} (2004), 585--609.

\bibitem{ISbea} A.\,L.\,Islas, C.\,M.\,Schober. Backward error analysis for multisymplectic discretizations of Hamiltonian PDEs. {\em Math. Comput. Simulations}, {\bf 69} (2005), 290--303.

\bibitem{KoideFur} S.\,Koide, D.\,Furihata. Nonlinear and Linear Conservative Finite Difference Schemes for Regularized Long Wave Equation. {\em Jpn. J. Ind. Appl. Math.}, {\bf 26} (2009), 15--40.

\bibitem{Kuperschmidt} B.\,A.\,Kuperschmidt. {\em Discrete Lax equations and differential-difference calculus}, Ast\'{e}risque No. 123, 1985.

\bibitem{LiSun} H.\,Li, J.\,Sun. A new multi-symplectic Euler box scheme for the BBM equation. {\em Math. Comp. Model.}, {\bf 58} (2013), 1489--1501.

\bibitem{MatsuoFur} T.\,Matsuo, D.\,Furihata. Dissipative or conservative finite-difference schemes for complex-valued nonlinear partial differential equations. {\em J. Comput. Phys.}, {\bf 171} (2001), 425--447.

\bibitem{McLQuisp} R.\,I.\,McLachlan, G.\,R.\,W.\,Quispel. Discrete gradient methods have an energy conservation law. {\em Discrete Contin. Dyn. Syst.}, {\bf 34} (2014), 1099--1104.

\bibitem{Robi} R.\,I.\,McLachlan, G.\,R.\,W.\,Quispel, N.\,Robidoux. Geometric integration using discrete gradients. {\em R. Soc. Lond. Philos. Trans. Ser. A Math.
Phys. Eng. Sci.}, {\bf 357} (1999), 1021--1045.

\bibitem{Olver} P.\,J.\,Olver. Euler operators and conservation laws of the BBM equation. {\em Math. Proc. Cambridge Philos. Soc.}, {\bf 85} (1979), 143--160.

\bibitem{Olverbook} P. J. Olver, {\it Applications of Lie Groups to Differential Equations, (2nd edn)}, New York: Springer Verlag, 1993.

\bibitem{Preissman}  A.\,Preissman. Propagation des intumescences dan les canaux et rivi\'{e}res. {\em First Congress French Association for Computation}. Grenoble, 1961.

\bibitem{S-SV86} J.\,M.\,Sanz-Serna, J.\,G.\,Verwer. Conservative and nonconservative schemes for the solution of the nonlinear Schr\"odinger equation. {\em IMA J. Numer. Anal.}, {\bf 6} (1986), 25--42.

\bibitem{Schober} C.\,M.\, Schober. Symplectic integrators for Ablowitz--Ladik discrete nonlinear Schr\"odinger equation. {\em Phys. Lett. A}, {\bf 259} (1999), 140--151

\bibitem{McLaren} G.\,R.\,W. Quispel, D.\,I.\, McLaren. A new class of energy-preserving numerical integration methods. {\em J. Phys. A}, {\bf 41} (2008), 045206.

\bibitem{QuispTurn} G.\,R.\,W.\,Quispel, G.\,S.\,Turner. Discrete gradient methods for solving ODEs numerically while preserving a first integral. {\em J. Phys. A.}, {\bf 29} (1996), 341--349.

\bibitem{SunQin} J.\,Q.\,Sun, M.\,Z.\,Qin. Multi-symplectic methods for the coupled 1D nonlinear Schr\"{o}dinger system. {\em Comput. Phys. Commun.}, {\bf 155} (2003), 221--235.

\bibitem{QinSun} Y.\,J.\,Sun, M.\,Z.\,Qin. A multi-symplectic scheme for RLW equation. {\em J. Comput. Math.}. {\bf 22}, (2004), 611--621.

\bibitem{Tang} Y.\,F.\,Tang, V.\,M.\,P\'{e}rez-Garc\'{i}a, L.\,V\'{a}zquez. Symplectic methods for the Ablowitz--Ladik model. {\em Appl. Math. Comput.}, {\bf 82} (1997), 17--38.

\bibitem{Vi1984} A.\,M.\,Vinogradov, The $\mathcal{C}$-spectral sequence, Lagrangian formalism and conservation laws I and II, {\it J. Math. Anal. Appl.}, {\bf 100} (1984), 1--129.


\end{thebibliography}
\end{document}